\documentclass[11pt]{article}
\usepackage{geometry}                
\usepackage{graphicx}
\usepackage{amssymb}
\usepackage{epstopdf}
\usepackage{psfrag}

\usepackage{amsthm}

\def\Bbb{\mathbb}
\def\R{{\Bbb R}}
\def\C{{\Bbb C}}
\def\Z{{\Bbb Z}}

\def\mbf0{\mathbf{0}}

\newcommand{\bmat}{\begin{pmatrix}}
\newcommand{\emat}{\end{pmatrix}}

\theoremstyle{theorem}
\newtheorem{thm}{Theorem}

\theoremstyle{theorem}

\theoremstyle{theorem}
\newtheorem{prop}{Proposition}

\theoremstyle{definition}

\theoremstyle{remark}

\title{The `Holiverse':    holistic   eversion of the 2-sphere in $\R^3$}

 \author{Iain R. Aitchison\\
 \\
 { Dedicated to Stephen Smale, who eschewed spin for deeper principles}}
 \date{}     

\begin{document}

\maketitle

\begin{quote}
{\small
\centerline{\bf Abstract}
We give a short, simple and conceptual proof, based on spin structures,  of sphere eversion:   
an embedded 2-sphere in $\R^3$ can be turned inside out by  regular homotopy. Ingredients of this eversion are seamlessly connected. We also give  the mathematical origins of the proof: the Hopf fibration, and the topological structure of real-projective 3-space.

}
\end{quote}

 \section{Introduction}
 
In 1957 Smale proved that any two smooth   immersions of the 2-sphere in $\R^3$ are regularly homotopic. As a corollary, any smoothly embedded sphere can be smoothly turned inside out by regular homotopy. Smale's proof is abstract, and far-reaching,  involving liftings and projections of homotopy fibrations, and for 50 years has inspired a number of attempts to conceptualize directly how such an eversion can be achieved, as well as several beautiful  computer animations of actual eversions. These animations have also provided the bench mark in computer-graphic realization of mathematical concepts.

We present here the first truly holistic eversion, by which we mean that each stage of the eversion is conceptually a natural part of a seamless whole. The origins of the eversion lie in the simple interplay of the Hopf fibration and the antipodal map on the 3-sphere, and the fact that an embedded  2-sphere in real projective 3-space $\R P^3$ can be everted essentially trivially. However, the actual eversion in $\R^3$
does not require a conceptual understanding of these origins to be understood in its entirety, and has its own intrinsic integrity. Accordingly we describe the actual eversion, and its mathematical origins, in two independent sections.

To provide conceptual context, we include a very brief historical account, mentioning  several of the ingredients of previous proofs. Essentially all of these proofs critically depend, at some stage, on direct visualization of part of the process, and in this regard differ from the proof offered here. Although visualization plays an important role in the communication of the essential ideas, this holistic eversion can be grasped conceptually, and as such is far  less reliant on the need for pictures or computer animation.
\vfill
\medskip
\hrule

\medskip
{\footnotesize
2010 Mathematics Subject Classification. Primary: 
57R42 ; 
secondary:  57M60, 
00A66.  
}

\pagebreak

 Essentially all existing animations rely at some stage in explaining the local topological changes occurring: for the holiverse, 
these details are unnecessary.
Moreover, the mathematical prerequisites for completely  understanding this holiversion are generally taught to students in advanced undergraduate or early graduate level topology classes. In spirit, this proof is closest to the first specific eversion described by Arnold Shapiro
\cite{FM}, although much simpler. We thank George Francis for suggesting the name `holiverse'.

\medskip

 \noindent{
 \bf A brief outline of the holiversion}: \  The eversion   arises from a  simple immersed disc, illustrated in Figure  \ref{fig:disctorusrot}: 
\begin{enumerate}

\item  {\bf The Disc}:  (a) Take a planar embedded disc, with two small distinguished arcs on the boundary: push these towards each other, twisting slightly while doing so --  as if to construct an annulus (0-twist) or Mobius band
($\pi$ twist). The boundary circle is drawn: you can imagine the disc conveniently wrapping around inside  a torus.
\subitem  (b) Push the arcs through each other slightly 
to  create the simplest immersed disc,  with double-point set a single arc. The projection of the disc to the plane is still an immersion:
 the boundary of this disc is still unknotted.  However, a collar neighbourhhod of the boundary circle   has two full negative twists
 when its core circle is stretched  out as a round planar circle.
\subitem (c) Smoothly unwist to   remove the innermost extraneous crossing, revealing  the boundary circle as a $(2,-1)$-torus knot; the collar annulus in the disc twists exactly as an annular neighbourhood of the circle on the torus. 
\subitem  (d) This introduces a `curtain-like'  bend in the disc, which can be done smoothly -- like the trace of the homotopy  $(x, x^3 + tx, t )$       viewed from the side $x$ axis. (The disc can be adjusted to meet the torus from outside, and  is slightly trickier to visualize.)   \begin{figure}[htbp] 
    \centering
 (1.1)  \includegraphics[width=2in]{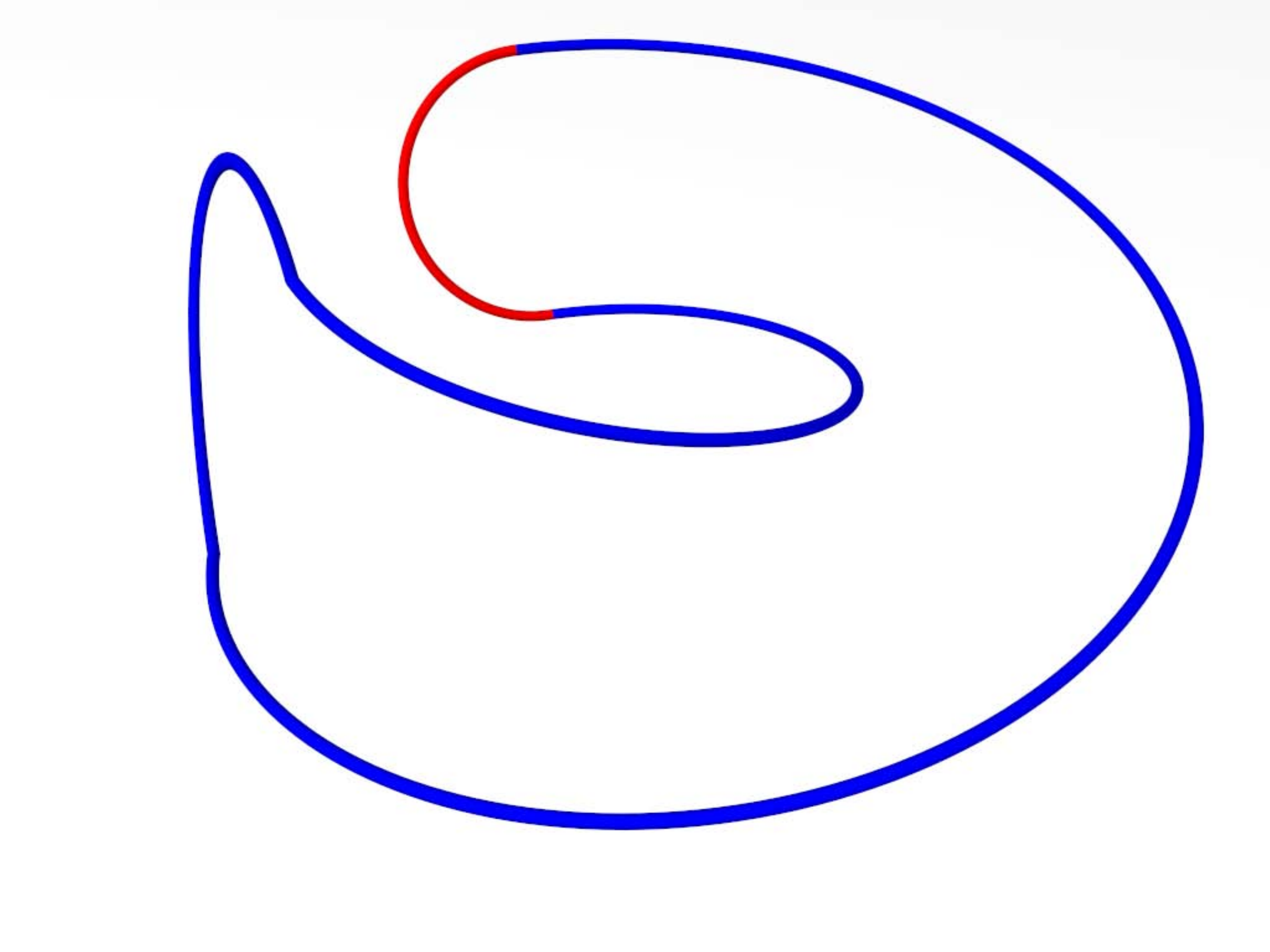}\qquad \qquad    \qquad
 (1.2) \includegraphics[width=2in]{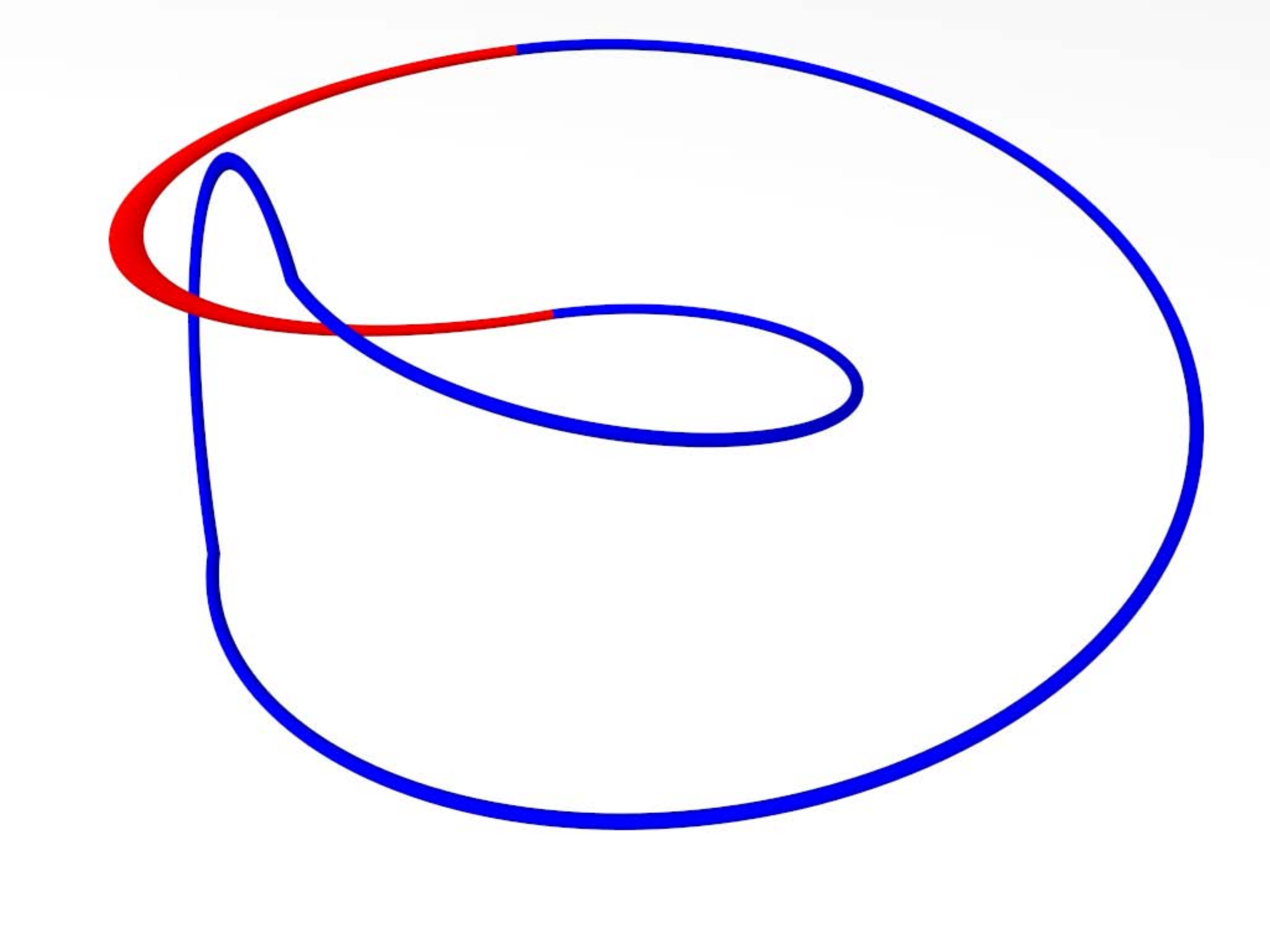} \hfill 
 \qquad  \qquad  \hskip 5cm
  (1.3) \includegraphics[width=2in]{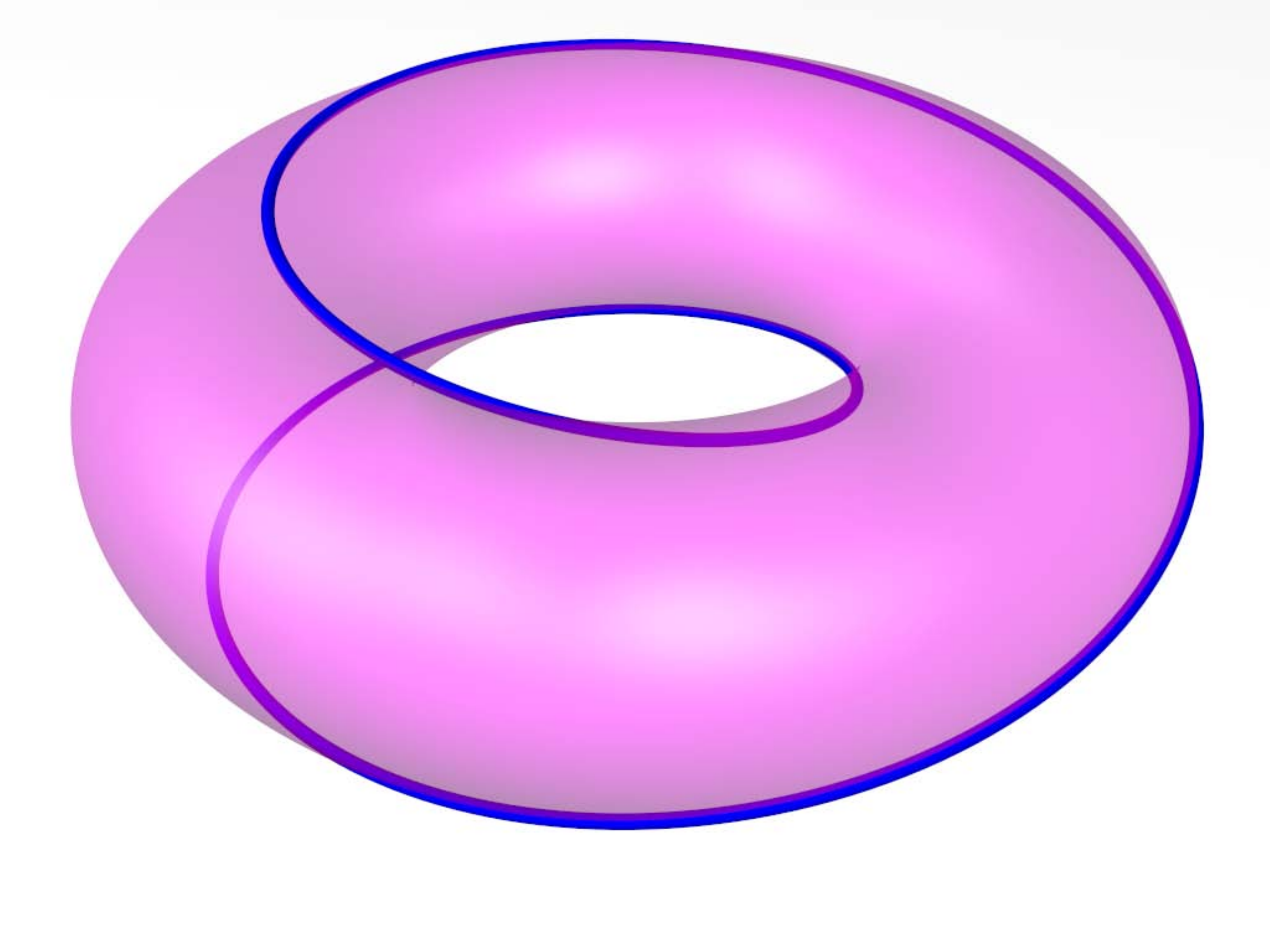} \qquad\qquad \qquad
(1.4)  \includegraphics[width=2in]{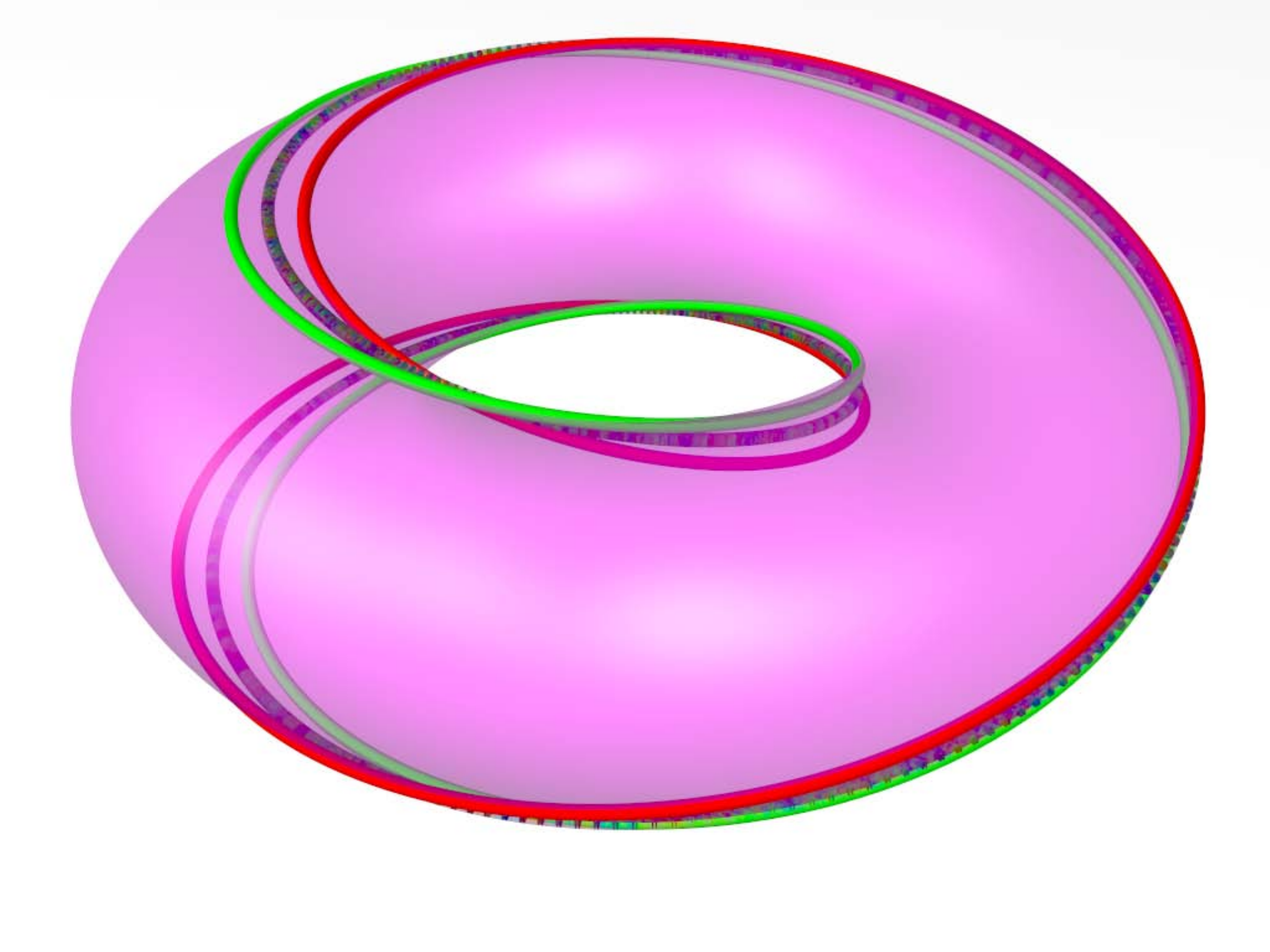}   \hfill\hskip 5cm
\caption{ (1.1)    An embedded disc  
 (1.2) An immersed disc, unknotted boundary circle. (1.3)  Untwisting gives   a $(2,-1)$ torus knot 
(1.4) Rotating two discs rigidly  in opposite directions.}
   \label{fig:disctorusrot}
   \end{figure}
    
\item {\bf  The Immersed Sphere:} (a) Thicken the disc to create an immersed sphere, which is essentially a union of two parallel copies of the disc, and an annulus following along the boundary circle, twisting to follow a collar in the disc orthogonally, with interior inside the torus.
\subitem (b) The sphere is regularly homotopic to an embedded sphere (shrink the disc).
\subitem (c) The discs and annulus have boundary   two parallel $(2,-1)$ circles.

\item   {\bf  The Eversion:} (a) Spin each disc rigidly in opposite directions  around the torus, so that  they coincide as they pass through  each other, and then return interchanged.
 
\subitem  (b) The two boundary circles bound an annulus, which sweeps across the torus, double covering a M\"obius band as the discs coincide, and returning to the initial confiuration, inside out with circles interchanged (and with a Dehn twist).

\subitem  (c) The two discs and annulus remain smoothly glued together: thus spinning turns the original immersed neighbourhood sphere inside out. Shrink to an embedding, to obtain an inside-out version of 2(b) above.

\end{enumerate}

 \section{Historical context:  underlying constructs}

We present a chronological synopsis of previous work on sphere eversion: this draws heavily on information organized by John Sullivan,
who credits George Francis as primary source \cite{Su}, and on correspondence with George Francis.

\begin{itemize}
\item 1924
J. W. Alexander \cite{Al} proves that every embedded 2-sphere in $\R^3$ is isotopic to the standard round sphere.

\item 1957: Stephen Smale announced his classification of immersions of the 2-sphere in $\R^3$ up to regular homotopy, 
published in \cite{Sm}, and generalizing the Whitney-Graustein classification of immersions of the circle in the plane. Smale's work   involves normalizing  a disc, homotopy fibrations, and homotopy classes of the tangent mapping into Stiefel and Grassman manifolds. 
Smale's approach immediately generalizes to higher dimension (Hirsch-Smale), and for immersions of arbitrary smooth manifolds. The original case
 implicitly   exploits $\pi_2(G) = 0$ for any Lie group $G$, in particular when $G= SO(3)\cong \R P^3 \cong L(2,-1)$, the Lens space. These spaces play a role underlying the holiverse: the author learned immersion theory from Smale in a graduate course at Berkeley in the early 1980's.

\item 1960. Arnold Shapiro described, but did not publish,  an  eversion using a neighbourhood of  Boy's immersion of the  real projective plane 
$\R P^2$ into $\R^3$, requiring a non-intuitive regularly homotopy, in several stages,  of an embedded sphere to the boundary of a twisted $I$-bundle neighbourhood of the immersed $\R P^2$. The projective plane is the union of a disc and a M\"obius band, a decomposition which features in the holiverse.
An exposition was given in 1979 by Francis and Morin \cite{FM}: the holiverse has several features in common, such as the appearance of the $(2,-1)$ torus knot, the need for an immersed disc bounding this knot, and the fact that a neighbourhood of an immersed disc immediately gives an immersed sphere regularly homotopic to an embedded sphere.
Nonetheless, Shapiro's eversion remained hard to visualize, especially so at the time it was announced: there were no computer aids for graphical representation.

\item 1966. Tony Phillips \cite{Ph} attempted to describe Shapiro's eversion pictorially, and in so doing obtained a new purely visual description of an eversion.
 Phillips' rendering of his procedure exploits sequences of pictures of immersed ribbons/annuli, in part with stages emulating the more familiar   Whitney-Graustein regular homotopies of circles in the plane.

 Around this time, in interaction with Froissart, 
Bernard Morin utilized the symmetry 
of a symmetric 4-lobed sphere immersion as half-way surface, as an alternative to an immersed projective plane; this makes manifest the equivalence of inside and out, with the  symmetry  interchanging these.   Again, an explicit non-intuitive regular homotopy from an embedded sphere must be demonstrated.
Fran\c cois  Ap\'ery and Morin subsequently showed that Morin's eversion has 
 the minimal number of topological events \cite{MP}.
 
\item 1968 Bryce De Witt \cite{BD} gave an outline of another pictorial scheme for an eversion, although details making this explicit have not appeared.

\item c. 1970 Charles Pugh, Smale's colleague  at Berkeley, constructed by hand a sequence of chicken-wire models showing the stages of Morin's eversion, based in turn on Morin's clay models. These were proudly suspended in the Mathematics Department until stolen sometime in the mid-late 1970's. 

\item 1974.  William Thurston, Hirsch's  student at Berkeley,  conceptualized an  eversion using the so-called
`belt trick', with origins in $Spin(3) = SU(2)$): a long strip can be given 
a full twist by either twisting the ends, or interchanging them by parallel
transport along a straight line. The resulting eversions use ideas of
 corrugations and symmetry, again with more general implications, but there also remains some  difficulty in  explicitly describing or following every stage of the eversion.

\item  1977.    Nelson Max spent six years digitalizing Charles Pugh's chicken-wire models, with coordinates for points on the models 
calculated by hand and entered as data for his landmark animated movie \cite{Ma}. This was an international sensation, and a tour de  force  in the dawning age of computer graphics: it motivated the author's interest.

\item 1979.  George Francis and Bernard Morin \cite{FM} published a description of  Arnold Shapiro's eversion of the sphere, pointing out the desire for a holistic eversion  in a similar spirit, with each stage understandable in terms of familiar topological concepts.

 \item 1987. George Francis'  `A Topological Picturebook' \cite{Fr} is published,
 with many hand-drawn pictures of stages of eversions.

\item 1992.  Fran\c ois Ap\'ery, with  Bernard Morin, describe an algebraic halfway model for the eversion of the sphere \cite{AM}. Being able to describe surfaces in such a way facilitates the production of computer graphics  and software such as Povray which render mathematical equations, rather than objects constructed from points in a large database.

\item 1994.    Silvio Levy, Delle Maxwell, and Tamara Munzner create the animation `Outside In', describing and explaining Thurston's ideas
\cite{OI}.

\item 1995.      John Sullivan, Bob Kusner and George Francis implemented a minimax eversion based on Kusner and Bryant's work on Willmore energies: Kusner, another Berkeley student, had found an analytic  surface which reminded him of Morin's half-way surface.
The use of a natural energy minimizing `Willmore flow' as a guiding principle recaptures earlier approaches more  beautifully, providing a heuristic justifying Morin's initial viewpoint. The Willmore flow clearly converges to an embedded sphere when run on Brakke's `Surface Evolver' program \cite{Bra}, but only recently have results on existence and uniqueness for the corresponding 4th-order PDE flow equations 
begun to emerge \cite{Su}.

\item 1998.     `The Optiverse',  an animation of the Willmore flow eversion created by 
     Sullivan,  Francis, and Stuart Levy, won a prize for animation at the 1998 Berlin International Congress of Mathematicians.  \cite{ FSH,FSK,FS, SFL} 
 \end{itemize}

\section{Holisitic  eversion in $\R^3$: Explicit   details }

Let  $D_{zx}$ be  the unit disc in the $zx$-plane with equation
$(x-\sqrt 2)^2 +z^2 = 1$, with  $S^1_{zx}  := \partial D_{zx}$  meeting the $+x$-axis at points $d_1<d_2$. Denote angle measure from the centre of $D_{zx}$  by $\phi$, measured anticlockwise from the $+x$-axis.
Let $SD_{z}\cong S^1\times D^2$ denote   the solid 
donut obtained by revolving $D_{zx}$ around the $z$-axis in $R^3_{xyz}$, while rotating the disc about its centre uniformly by $\pi$ anticlockwise, so that $\{ d_1,d_2\}$ creates a $(2,-1)$ torus knot $K_0$. Foliate $D_{zx}-\{ d_1,d_2\}$ by arcs of circles whose centres lie on the $x$-axis, and which are orthogonal to $S^1_{zx} $. Let $\alpha_\phi$ denote such an   arc with endpoints on $S^1_{zx}$ at angles $\pm \phi$, $0<\phi <\pi$, and let   $A_\phi\subset SD_{zx}$ denote the annulus obtained by rotating   $\{\alpha_{\phi},\alpha_{\pi -\phi} \}$; for $ \phi =   \pi/2$, this covers a M\"obius band $M$. Let $\theta$ denote standard angle measure in the $xy$-plane, and let $K_\theta$ denote the image of $K_0$ after $\theta$-rotation about the $z$-axis, $-\pi\leq \theta \leq \pi$. 
 $K_\theta$ meets $S^1_{zx} $ in two antipodal points $ S^0_\theta$, at angles $\phi = -\theta/2, \pi-\theta/2$; $K_\theta\cup K_{-\theta}  = \partial  A_{\theta/2} \subset  \partial SD_{z}$. 
 
\begin{prop}
There exists a smoothly immersed disc $D_0\subset R^3_{xyz}$,  $\partial D_0 = K_0$,    orthogonal to $\partial SD_{z}$, meeting from the outside.
\end{prop}

\noindent
{Proof.}   Take an equilateral triangle in $R^2_{xy}$, with one vertex on the $+x$-axis and one edge on the $y$-axis. Fold over a square 
attached to the top edge, on the $+z$ side; fold under another square 
attached to the bottom edge, on the $-z$ side, to create an embedded disc which projects to the plane with two fold curves. Change the crossing between the two edges of folded squares
which project to intersect on the $-x$-axis. Smooth the construction 
to obtain an immersed disc $D'$ with unknotted boundary and a proper  arc of intersection, and observe that each of (i) a collar neighbourhood of the boundary circle of $D'$, (ii) a collar neighbourhood of $K_{\pm \pi}$ in $M$, and (iii) a neighbourhood of $K_0$ on $\partial SD_{z}$, is an unknotted annulus with four negative half-twists. Using a collar of $\partial SD_{z}$ and a regular homotopy/isotopy of $D'$, we obtain the desired immersed disc $D_0$, with boundary-circle-collar embedded in $R^3_{xyz}-int(SD_z)$.
Let $D_\theta$ be the $\theta$-rotation of $D_0$, with $\partial D_\theta = K_\theta$, and define $S_\theta := D_{-\theta} \cup A_{\theta/2} \cup D_{\theta}.$ 

\begin{prop}
$S_\theta$ is regularly homotopic to an embedded sphere in $R^3$, for (non-zero) $  \theta \in (-\pi,  \pi)$.  
\end{prop}

\noindent
{Proof.}  Let $S^1_t$ be the circle of radius $t$ in the $uv$-plane,   $\, 0< \epsilon \leq t \leq 1+\epsilon$, bounding the disc $D^2_t$, so that $D^2_{1+\epsilon}$ gives a neighbourhood of the unit disc $D^2_1$. 
The 2-spheres
$S^2_t= \partial (D^2_t\times [-t,t]) := \partial N_t \subset R^3_{uvw}$ are concentric with the boundary $S^2_{1+\epsilon}$ of a tubular neighbourhood of $D^2_1$ in $R^3_{uvw}$. These spheres have a natural splitting as two discs $D^\pm_t$ and an annulus $\bar A_t$: the construction can naturally be done smoothly. Doing so, let $\rho : D^2_1 \to R^3_{xyz}$ be a smooth immersion with image $D_0$, 
and $\rho^* : B^3\cong N_{1+\epsilon}=D^2_1\times [-1-\epsilon, 1+\epsilon] \to R^3_{xyz}$ be an immersion with image an immersed neighbourhood of the immersed disc $D_0$, with $\rho^* (\bar A_{1+\epsilon} )= A_{\theta_0/2}$ for some small $\theta_0$. We may assume that $\rho^* (S^2_\epsilon) $ is smoothly embedded; $\rho^*(S^2_t)$
gives a regular homotopy to $\rho^* (S^2_{1+\epsilon} )
:= D^-_ {1+\epsilon} \cup A_{\theta_0/2} \cup D^+_ {1+\epsilon} $.
Both $D^{\pm}_{1+\epsilon}$ are regularly homotopic to $D_0$, and 
have boundary $K_{\pm \theta_0}$;  they are regularly homotopic keeping $K_{\pm \theta_0}$ fixed, to $D_{\pm \theta_0}$. 
Hence 
$\rho^* (S^2_{1+\epsilon} )$ is regularly homotopic to $S^2_{\theta_0}$. 

\begin{prop}
$S_{\theta_0}$ is regularly homotopic to itself with orientation reversed.    
\end{prop}

\noindent
{Proof.} The regular homotopy is provided by the sequence of spheres $S_{\theta},\ \theta \in [\theta_0, 2\pi  - \theta_0]$. The intermediate sphere $S_{\pi}$ double covers the immersed projective plane 
$M\cup D_{ \pi}$. Doing this continuously, observe that the two discs $D_{\pm \theta_0}$ are interchanged after being rotated in opposite directions; the annulus $A_{\theta_0/2}$ undergoes an additional Dehn twist.

\begin{thm}
{\rm (Smale \cite{Sm})} 
An embedded sphere in $R^3$ is regularly homotopic to itself, reversing   orientation.
\end{thm}
 
\medskip

\noindent
{Proof.} This is an immediate application of the preceding propositions:
An embedded disc expands to create a single arc of self-intersection,
with the boundary of the immersed disc still unknotted but embedded on the torus as a $(2,-1)$ torus knot. An immersed 2-sphere neighbourhood of the disc is regularly homotopic to an embedded sphere, using neighbourhoods of immersed discs as they shrink back to the original disc. This immersed sphere is naturally the union of two `parallel' copies of the disc, and an equatorial annulus which is an embedded annulus in the solid torus, with boundary two parallel copies of the $(2,-1)$ torus knot. Rotating both of these circles in opposite directions, we can rotate the two immersed hemispheres of the disc to interchange. Simultaneously the annulus enlarges inside the solid donut, passes through itself as a double-covered M\"obius band, and then returns to the original annulus (after a Dehn twist). The sphere has now turned itself inside out smoothly, and can be shrunk back to an inside-out embedded sphere, completing the eversion.
There is an immersed projective plane at the half way stage. 
\bigskip

The theorem is a corollary of the more general results   first published in \cite{Sm}.

 \bigskip
 The mathematics underlying this proof simply combines basic facts about $Spin(3)$, the Hopf fibration, the Clifford torus, $SO(3)$,  $\R P^3$, the Lens space $L(2,-1)$, the simplest embedded M\"obius band and the simplest (non-trivial) immersed disc in $\R^3$, and $S^3$ as $-1$ surgery on the unknot:   this is explained in the next section. 

   \section{Eversion by spin: the origin of the eversion in $\R P^3$}
  
 A 2-sphere can be trivially turned inside out in $\R P^3 = B^3\cup N(\R P^2)$, the union of a ball and a twisted $I$-bundle neighbourhood of a 
 projective plane $\R P^2$:
radially expand an embedded sphere metrically centered at the centre of the ball until it becomes the boundary of $N(\R P^2)$. Push the sphere across the projective plane, so that it ceases to be an embedding only at the instant it double covers the projective plane. It then becomes the boundary of the ball, but with its orientation reversed. Since the 2-sphere is simply connected, this regular homotopy lifts to the universal (2-fold) cover $S^3$ of $\R P^3 = S^3/a$, with exactly two lifts since $\pi_1(\R P^3)$ is $\Z/2\Z$: we see two 2-spheres, centered respectively at $0$ and $\infty$ in $S^3\equiv \R^3\cup \infty$. These bound $S^2\times I$, and $\{ S_t \} := \{ S^2\times \{ t\} \}$ determines the  regular homotopy  to interchange the 2-spheres
$\{ S_{\pm \epsilon} \}  $ by isotopy, turning a single 2-sphere inside-out  in $\R P^3$. 

We modify the
radial expansion of $S^2$ by twisting upper and lower hemispheres as we radially expand, in a manner prescribed by the Hopf fibration. This will allow the regular homotopy in $\R P^3$ to be realized in $\R^3$, by first constructing a map 
from $\R P^3$ into $R^3$  
whose composition 
with the regular homotopy  
is a sphere eversion in $\R^3$. This map, naturally constructed using the Lens space structure of $\R P^3$, intertwines the Hopf flow  with the standard rotation of $R^3$ round the $z$-axis, making 
the complete eversion comprehensible and describable in `one shot'.  

$S^3 =\partial B^4\subset
\R^4\cong  \C^2$  
is defined by $\{ (\alpha,\beta)\, |\,  \alpha\bar\alpha +\beta\bar\beta=1\}$.  
   The subset $\alpha\bar\alpha = \beta\bar\beta$  is the 
   Clifford torus $CT \cong S^1_{1/2} \times S^1_{\scriptsize  1/2}$.
$CT$ defines   two complementary `solid donuts' $SD , \ SD^* \subset S^3$ with core circles  
$\{ (e^{i\phi},0)\}$, $\{ (0, e^{i\psi})\}$ in the $\alpha$ and $\beta$ complex lines,  
giving the standard genus 1 Heegaard decomposition of $S^3$. 
   Define oriented circles  
on $CT$  
   by   $ \ \lambda   
   :=\{ (e^{i \phi}, 1)/\sqrt 2\}, \mu   
   := \{ (1, e^{i \psi})/\sqrt 2 \},$   meeting at  
   $P\in CT$, and which give   homology basis $[\lambda]= (1,0),\  [\mu]=(0,1)$ respectively for $H_1(CT; \Z)$.
Thus $SD \cong \lambda \times D$, where $\partial D = \mu$, and 
$SD^* \cong   D^*  \times \mu$, with $\partial D^* = \lambda$, where $D,\ D^* $     
are 2-dimensional discs:  
 solid donuts  are parametrized by   $S^1$-actions on $S^3$ with   circle-of-fixed-points $\lambda$ or $\mu$. Identify 
$ e^{i \psi}$ with $ (1, e^{i \psi})/\sqrt 2 \in \mu \subset  \C^2$, corresponding to the $S^1$-action $(\alpha,\beta) \to e^{i\psi}(\alpha,\beta) :=  
(\alpha, e^{i\psi}\beta)$, and
define
$D^*_\psi := (   D^* ,  e^{i \psi} ) \subset SD^*$. Thus $P$ has orbit $\mu$, and $\partial D^*_\psi := K^*_\psi \subset CT$ is a circle isotopic to $ K^*_0 = \lambda$.
$SD^*$   constitutes a 2-handle and 3-handle attached to $SD$ to give $S^3$:   the 2-handle $\cup_{|\psi|\leq \epsilon} \, D^*_\psi$ has core disc $D^*_0$ attached along $\lambda$.

One lift of a radial regular homotopy of $S^2$ in $\R P^3$ begins with radial expansion of a small sphere $\bar S_*$ centered at $0 := D^*_0 \cap \{ \beta = 1 \}$, until it becomes $\bar S^2_0$, tangent to $CT$ along $K^*_0$. Let $\bar S^*_\psi$, $0<\psi<\pi$,
denote the expanded sphere intersecting $CT$ in the pair of circles $K^*_{\pm \psi}$: these circles also bound an annulus $\bar A^*_\psi$ in $SD$, which decomposes the 2-sphere into
$\bar  S^*_\psi = \bar D^2_\psi\cup \bar A^*_\psi \cup \bar
D^2_{-\psi}$.  The  equatorial sphere $\bar S^2_{\pi/2} \subset S^3$    lies midway in the isotopy of $\bar S_*$ to the antipodal small sphere $a( \bar S_*)$, double-covers a projective plane in $\R P^3$, and
  bounds a fundamental domain for $a:  S^3\to S^3$.  
  
The circles $K^*_{\pm \psi}$ also bound
the pair of discs $D^*_{\pm  \psi} \subset SD^*$,  meeting $CT$ orthogonally. In $SD$, choose   annuli $A^*_ \psi$, $ \psi \in (0,\pi)$,   meeting $CT$ orthogonally along $K^*_{\pm  \psi}$, intersecting each disc of $SD$ in a circular arc,  and define $S^*_ \psi := 
D^*_{ \psi} \cup A^*_ \psi \cup D^*_{- \psi}$:  each sphere
$S^*_{\psi} $ is the boundary of a 3-ball neighbourhood of the disc $D^*_0$, the core disc of the 2-handle, and there is a natural collapse of the disjoint family $\{ S^*_\psi \}$ onto $D^*_0$. This family also projects to give an eversion in $\R P^3$: We use the Hopf fibration to modify this eversion.  

   The Hopf flow and Hopf fibration is determined by the circle action  $e^{i\theta}_{\rm Hopf}(\alpha, \beta) := (e^{i\theta}\alpha,  e^{i\theta }\beta  )$.  
    Each of $SD, \ CT, $ and $SD^*$ is   invariant under the Hopf flow, and hence under the antipodal map $a$, which corresponds to $\theta = \pi$.

The Hopf circles on the Clifford torus are orbits of this free circle action, and in coordinates, are $(1,1)$ circles. 
   The   Hopf circle through $P$ intersects  $K^*_\psi$
   at  
$P^*_\psi = (e^{i\psi}, e^{i\psi}     )/\sqrt 2$. 
Let $D^h_\theta    $ denote the rotated disc $e^{i\theta}_{\rm Hopf}
(D^*_\theta)$: the point $P\in D^*$ follows $P^*_\theta$,   inducing a twisted reparametrization $A^h_\theta$ of  $A^*_\theta$ as we isotope  the spheres $S^*_\theta$ across $SD$.

Change   parametrization of $SD$ by a  
 left-hand Dehn twist, so that the Hopf flow on $SD\cong S^1\times D^2$ is 
$e^{i\theta}_{rot}(x,y) = (e^{i\theta }x, y)$. (This corresponds to constructing $S^3$ by $(-1)$-surgery on the unknot.) The $(1,1)$ Hopf circle through $P$  becomes a $(1,0)$ circle, $P^*_\theta$ lies on this circle,
and $\lambda$  becomes a  $(1,-1)$ circle, so that the discs 
$D^h_\theta \cong D^*_\theta    $ constituting $SD^*$ are attached as Hopf-rotated discs $e^{i\theta}_{rot}(D^h_0) = D^h_\theta$.  
The rotated discs $D^h_{\pm \pi/2}$ are equivariantly interchanged by $a$, and the annulus $A^h_{\pm\pi/2}$ double covers a M\"obius band in $\R P^3$. 
 Let 
$S^h_\theta =
D^h_{\theta}\cup A^h_{\theta}\cup  D^h_{-\theta} $ denote the    twist-parametrized spheres, and 
$S^2_\theta =D^2_{\theta}\cup A^2_{\theta/2}\cup  D^2_{-\theta} $ their projection  to $\R P^3 ={SD/a}\,\cup\, {SD^*/a}$. The eversion in $\R P^3$ can be constructed from a ball neighbourhood of $D^2_0$ and $\{ S^2_\theta\}$,
$\theta \in (0,\pi)$. Observe $A^2_{\theta/2}\subset {SD}/a$;
$D^2_\theta = e^{i\theta}_{rot} D^2_0$. 

The solid donut $SD/a$ also embeds in $\R^3$ realizing the Hopf circle action as rotation around the $z$-axis; the attaching circle for the 2-handle arising from $SD^*/a$ is a $(2,-1)$-torus knot $K_0$, which bounds a disc $D_0$ immersed in $\R^3$. Thus the discs $D_\theta$ are all embedded by rotating $D_0$, and the spheres $S^2_\theta$ all simultaneously smoothly immerse in $\R P^3$. Eversion in an abstract $\R P^3$  is emulated in the corresponding (non-immersed) image of $\R P^3$ in $\R^3$  described in the previous section.

\section{Additional explanatory figures}

   \begin{figure}[htbp]  
       \centering
 \includegraphics[width=2.8in]{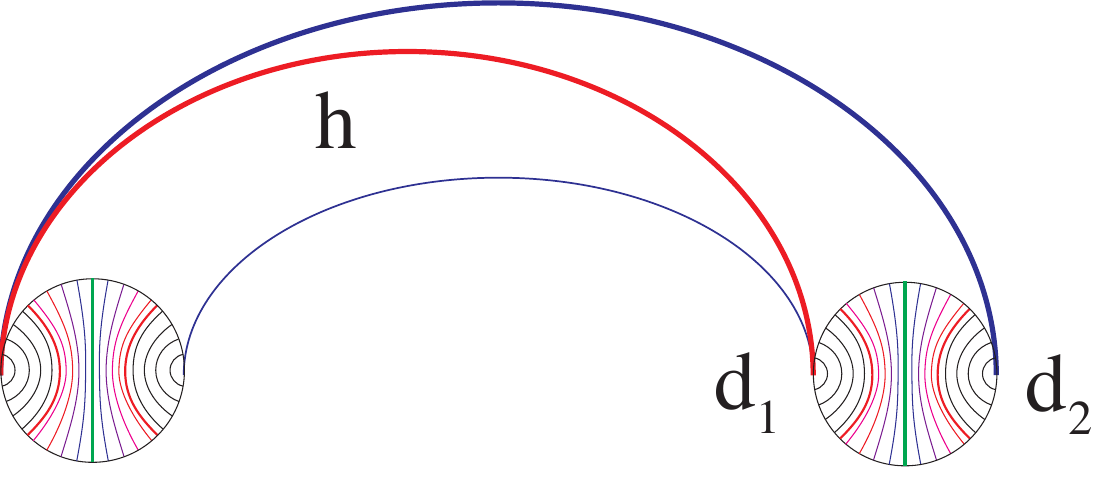} \  
  \includegraphics[width=2.8in]{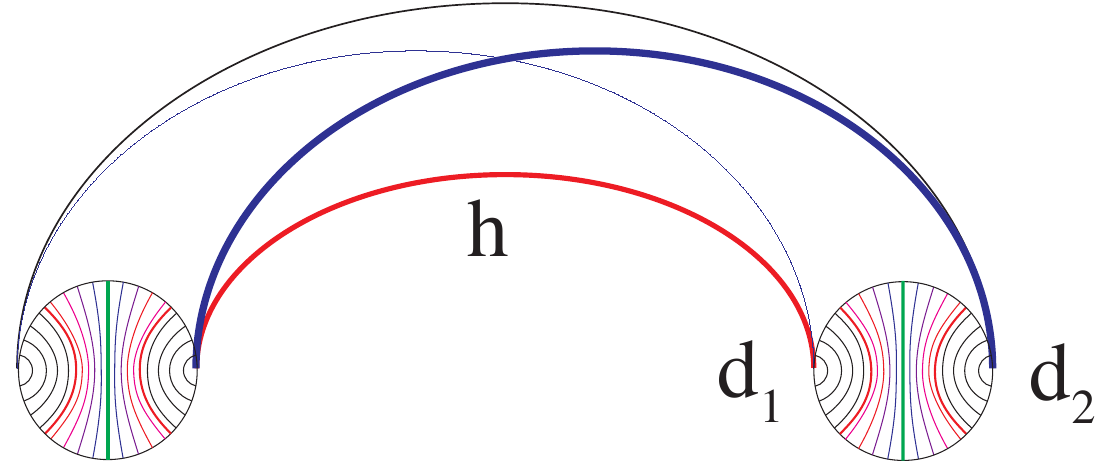} 
     \caption{ 
     The solid donut $SD_z$ obtained from $I\times D^2$ by identifying 
     two end-discs after $\pi$-rotation. The two points $d_1,d_2$ and arcs in the foliated disc sweep out the $(2,-1)$ knot $K_0$, annuli $A_\phi$ wrapping twice, and a M\"obius band $M$ created from the vertical arc.
($h$ stands for `Hopf circle' - see later)    }
    \label{fig:1121-1byes}
 \end{figure}
 
  \begin{figure}[htbp] 
    \centering
 \includegraphics[width=1.7in]{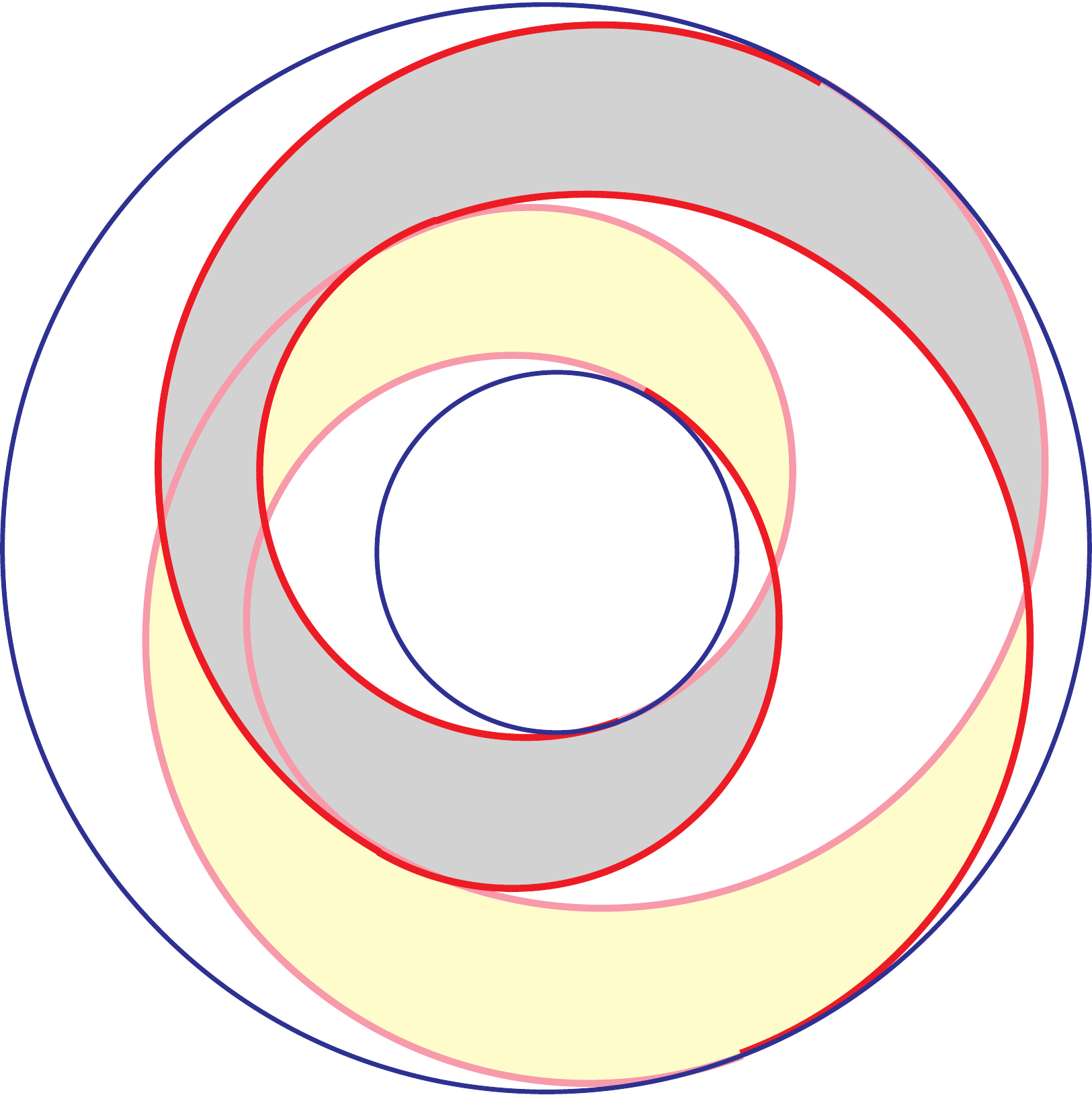} 
 \qquad\qquad
 \includegraphics[width=1.7in]{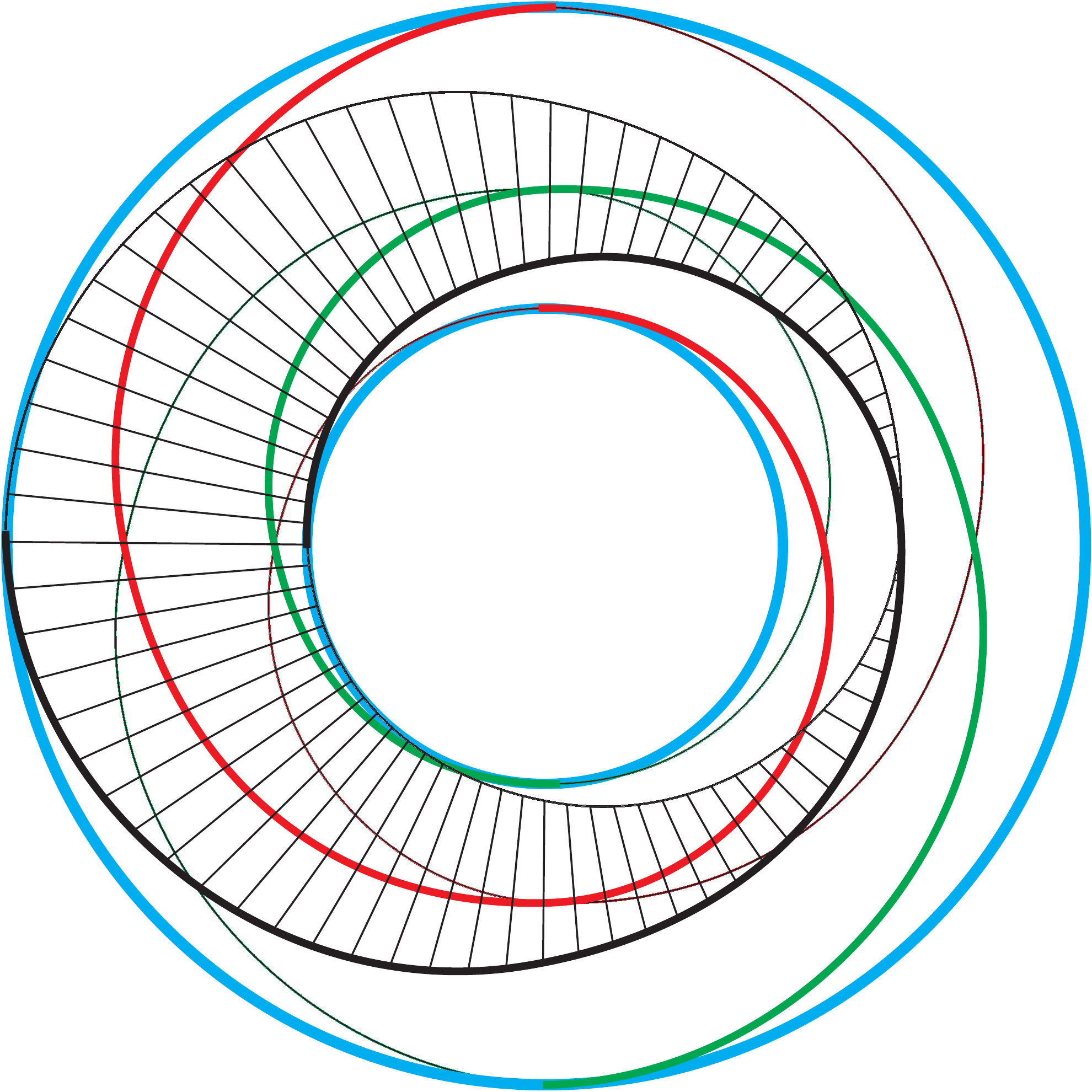} 
     \caption{ 
 An annulus $A_\phi$ on the left; on the right,  two knots $K_{\pm \pi/2}$ and $M$, with $\partial M = K_{+ \pi}\equiv K_{-\pi}$.  }
    \label{fig:ever3}
 \end{figure}
 
 \begin{figure}[htbp]  
    \centering
 \includegraphics[width=6in]{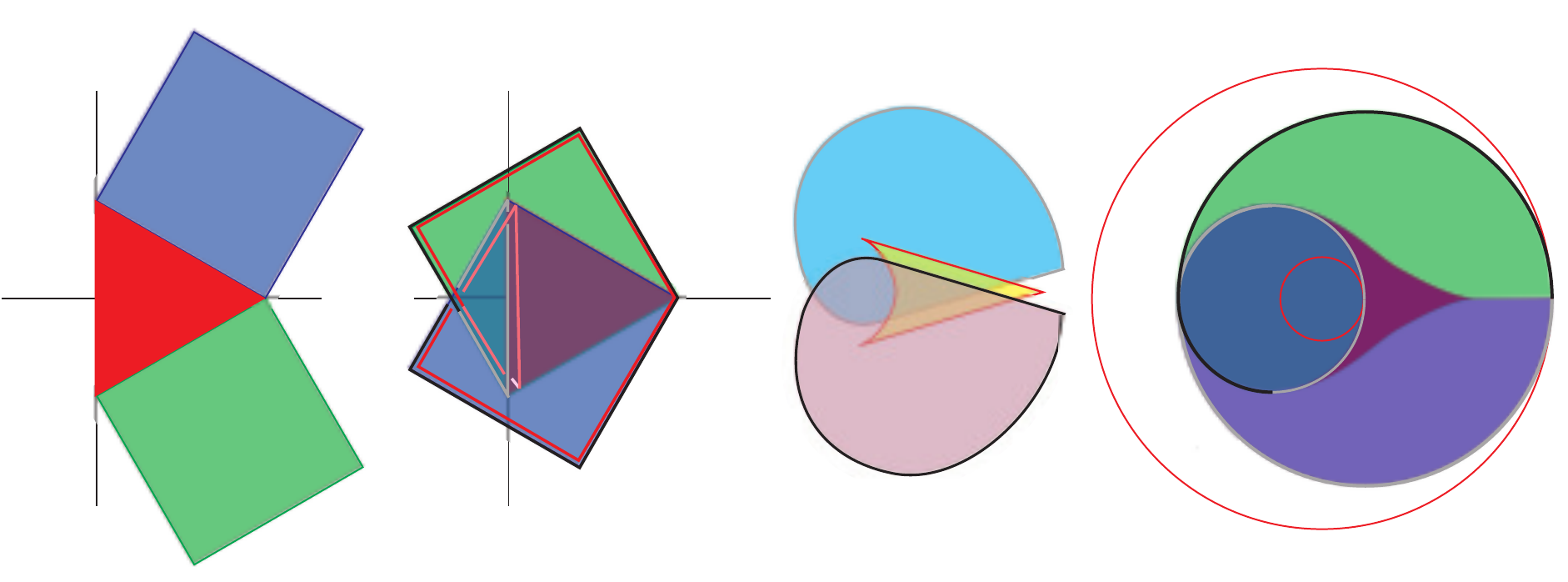} 
     \caption{ 
   Folding two squares over a triangle, and intersecting, before smoothing   
   to create $D'$ and $D_0$. Alternatively, create an intersection first, and then `fold', as in Figure \ref{fig:disctorusrot}. }
    \label{fig:SquareSmooth}
 \end{figure}
 
 \begin{figure}[htbp]  
    \centering
 \includegraphics[width= 5in]{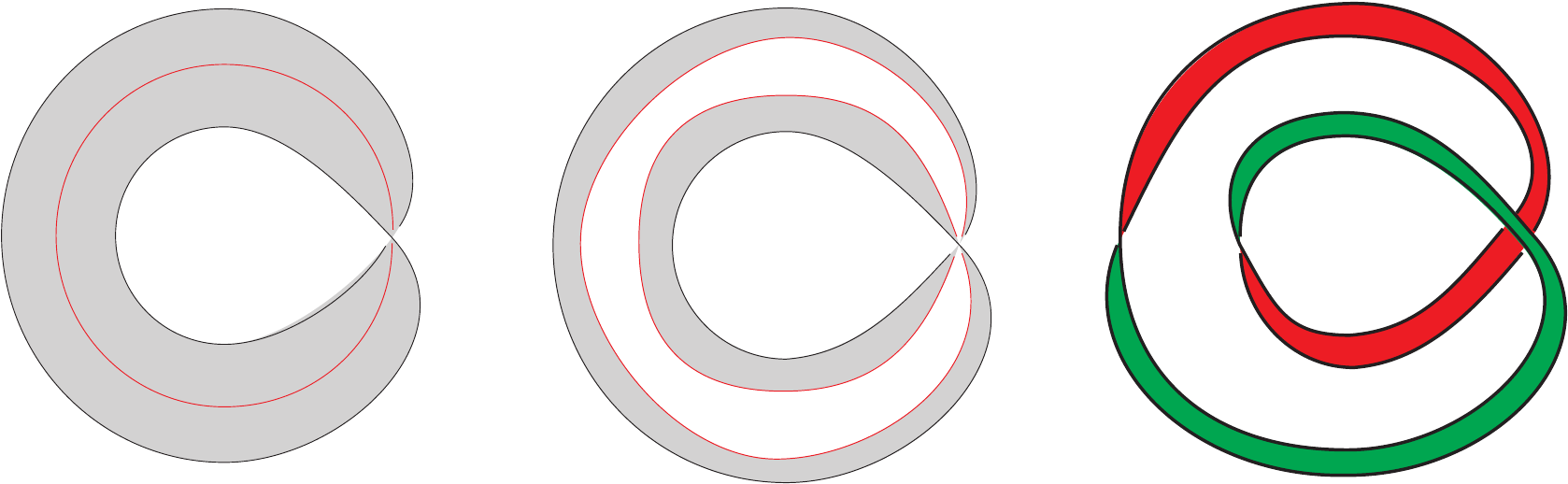} 
     \caption{ 
 The $(-4)$-half-twisted annulus as the collar of the M\"obius band, 
 also appearing in Figure \ref{fig:SquareSmooth} as the collar of the immersed disc, and as an annular neighbouthood of the $(2,-1)$-knot on the torus.  }
    \label{fig:Moanulusa}
 \end{figure}

   \begin{figure}[htbp] 
       \centering
 (a)   \includegraphics[width=2.2in]{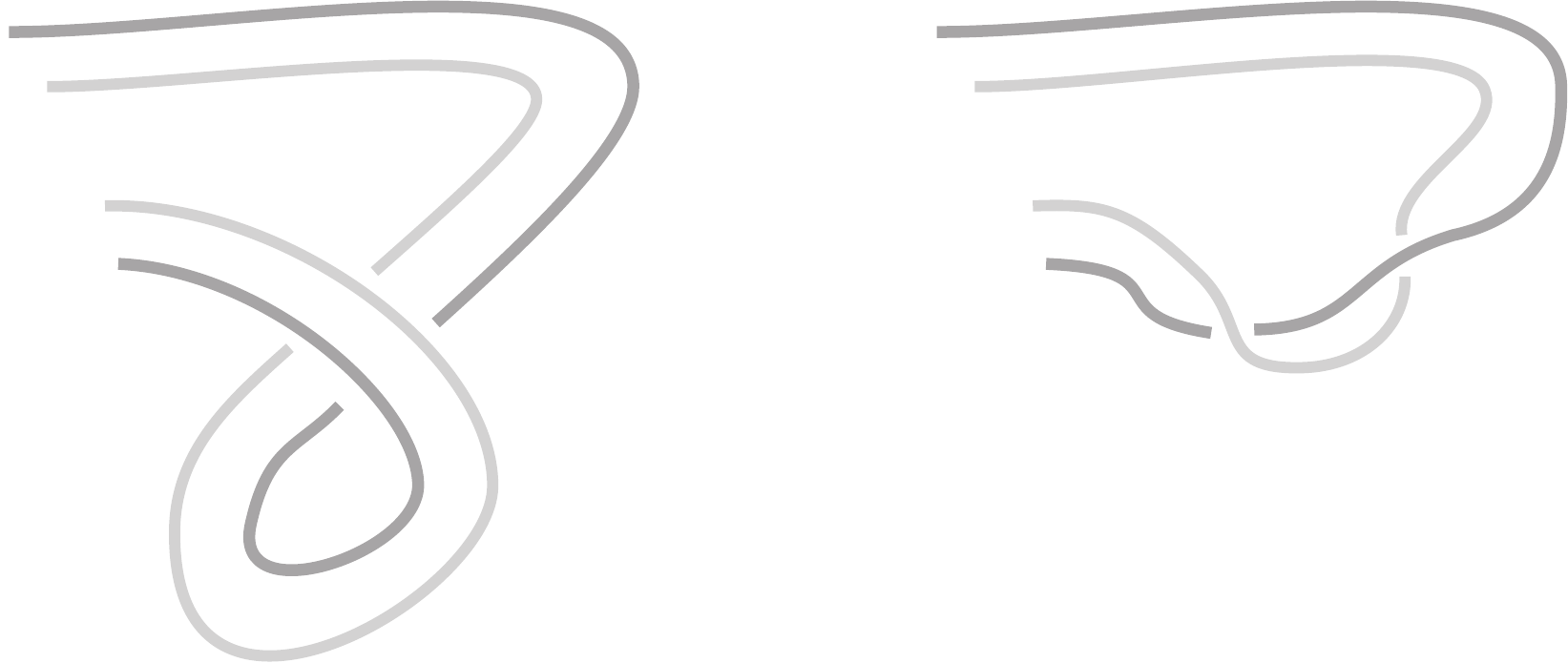} 
\qquad\qquad 
 (b)   \includegraphics[width=2.2in]{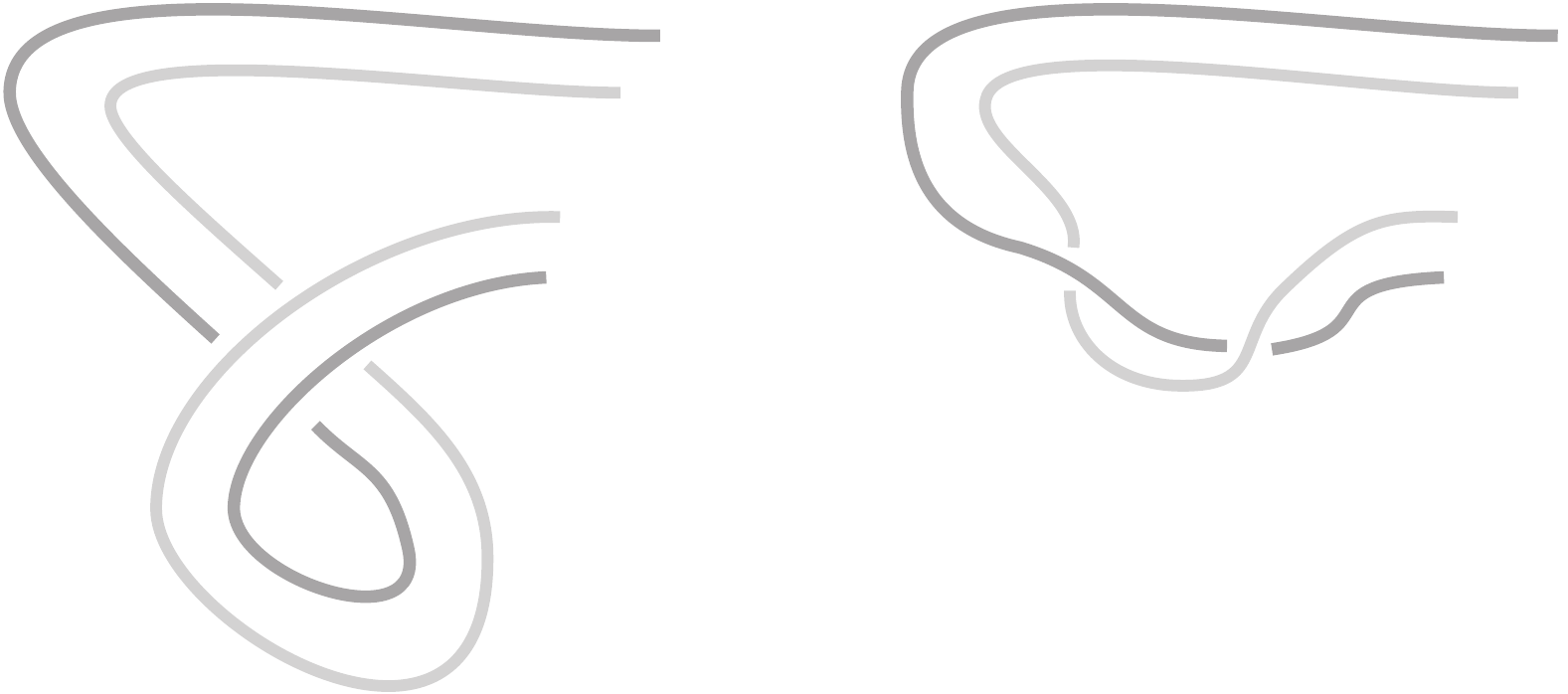} 
    \caption{ 
    A band in $\R^3$ with loops can always be pulled tight, converting loops into an even number of half-twists. (b) A loop with a negative crossing yields two negative half-twists. Positive and negative loops are regularly homotopic.}
    \label{fig:minustwista}
    \centering
(i) \includegraphics[width=1in]{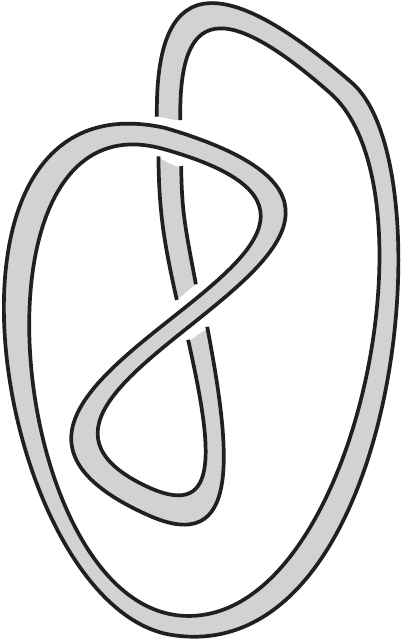} \qquad  
(ii)  \includegraphics[width=1in]{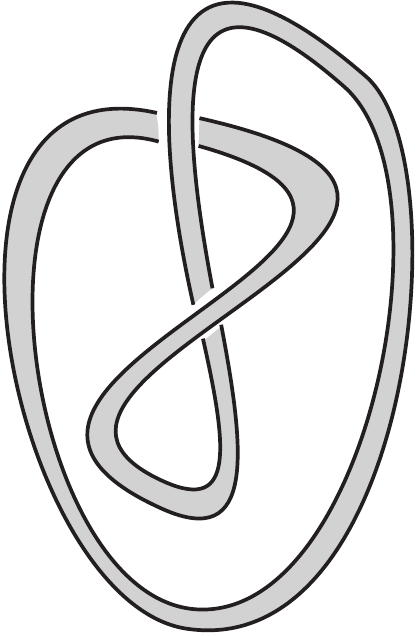}  \qquad 
(iii)   \includegraphics[width=1in]{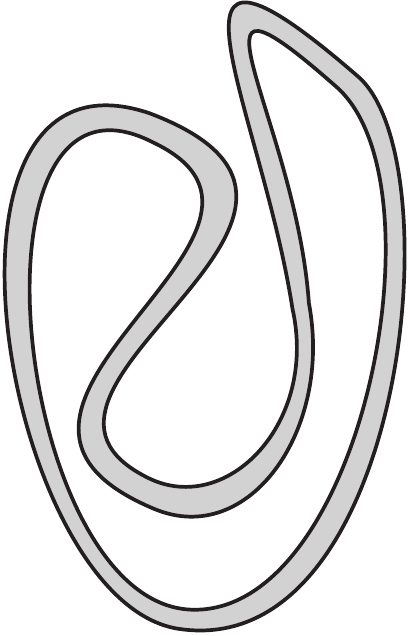} 
     \caption{ 
Untwisting the (-4)-half-twisted annulus in  (ii) by regular homotopy in $\R^3$, maintaining planarity of projection and creating a disc with a single clasp intersection:
(i) vertical  crossing change  homotopy (iii) horizontal homotopy }
    \label{fig:Untwist}
    \centering
 (a)   \includegraphics[width=1in]{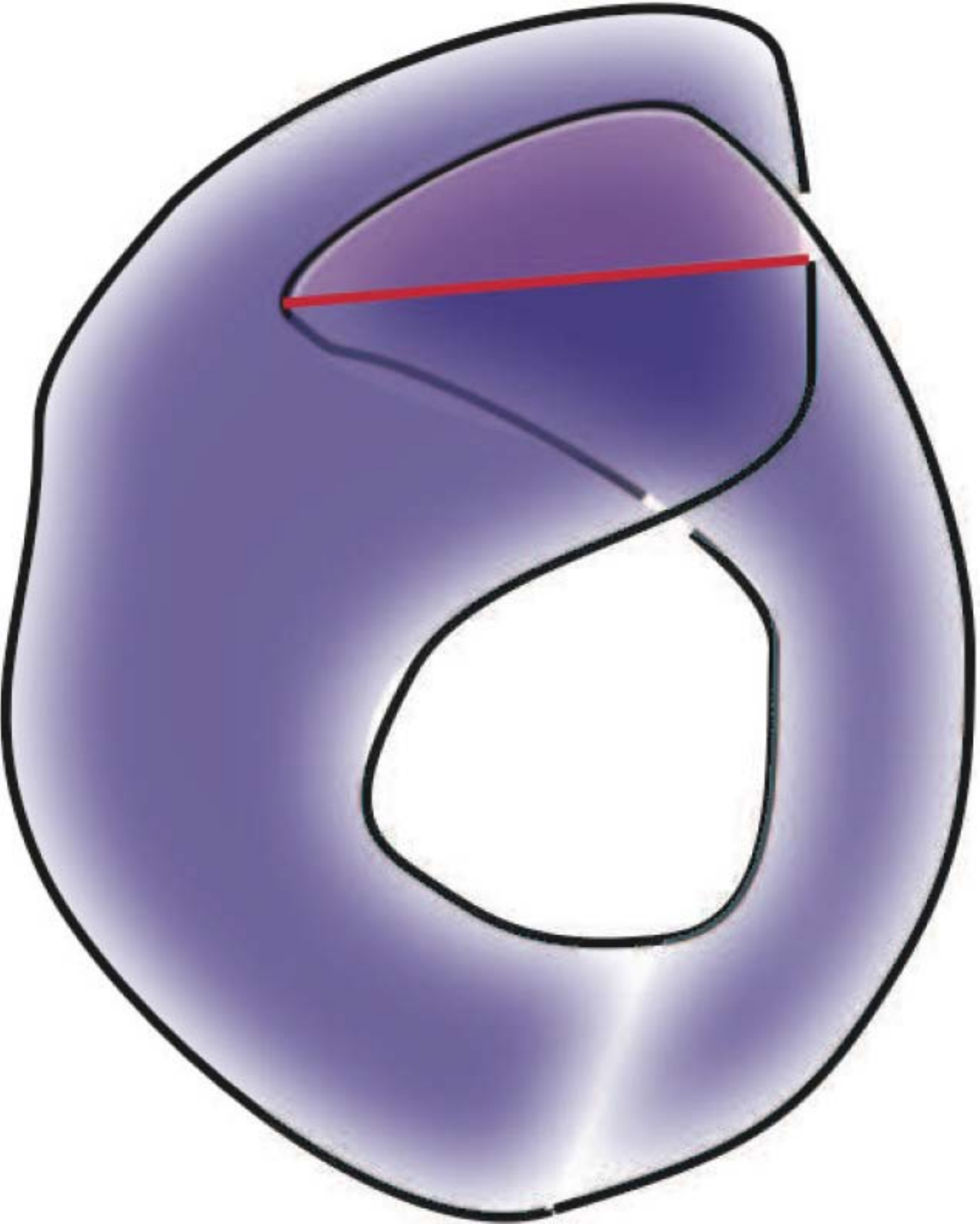}  
 (b)   \includegraphics[width=1in]{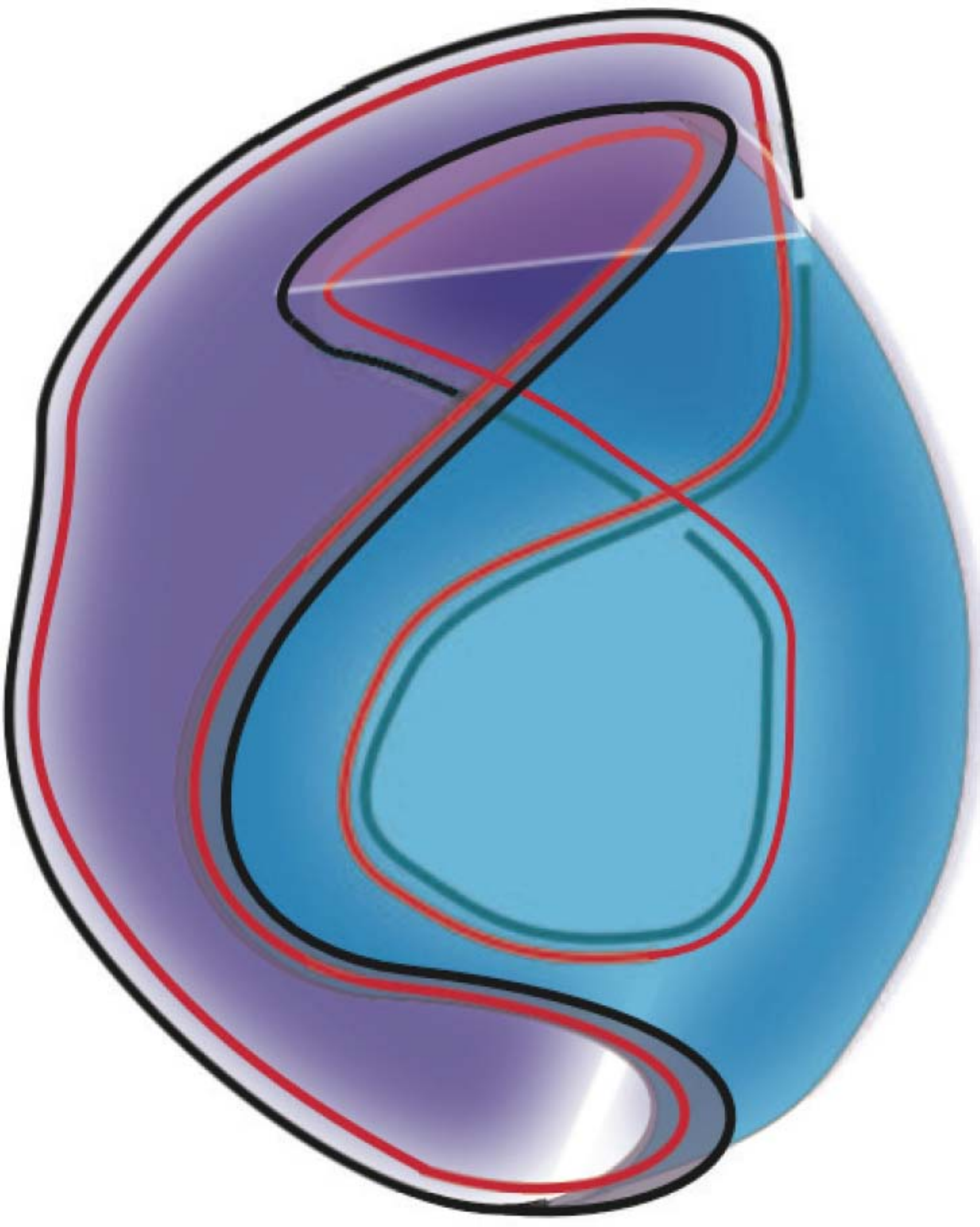}
  (c)   \includegraphics[width=1in]{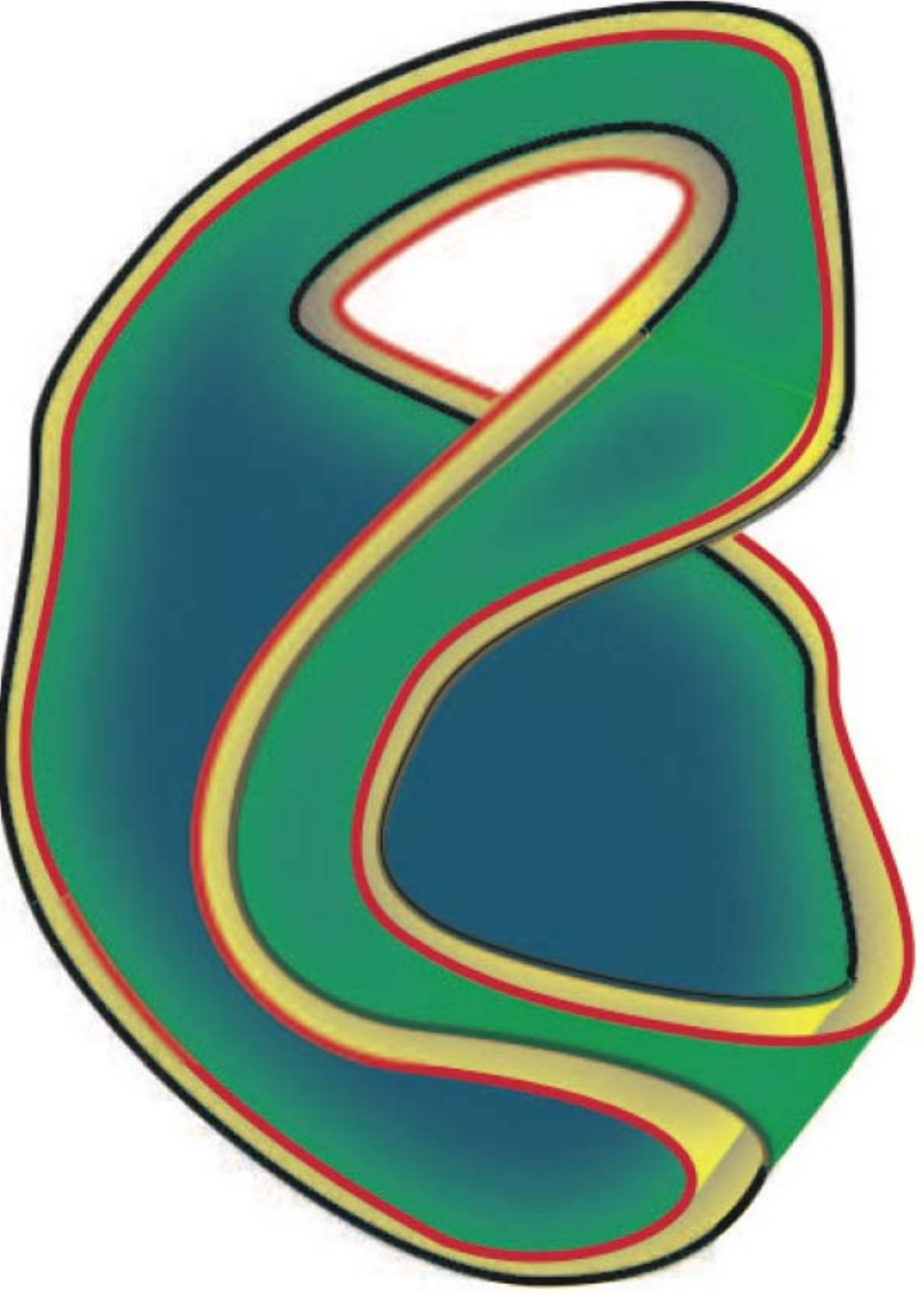} 
(d)    \includegraphics[width=1in]{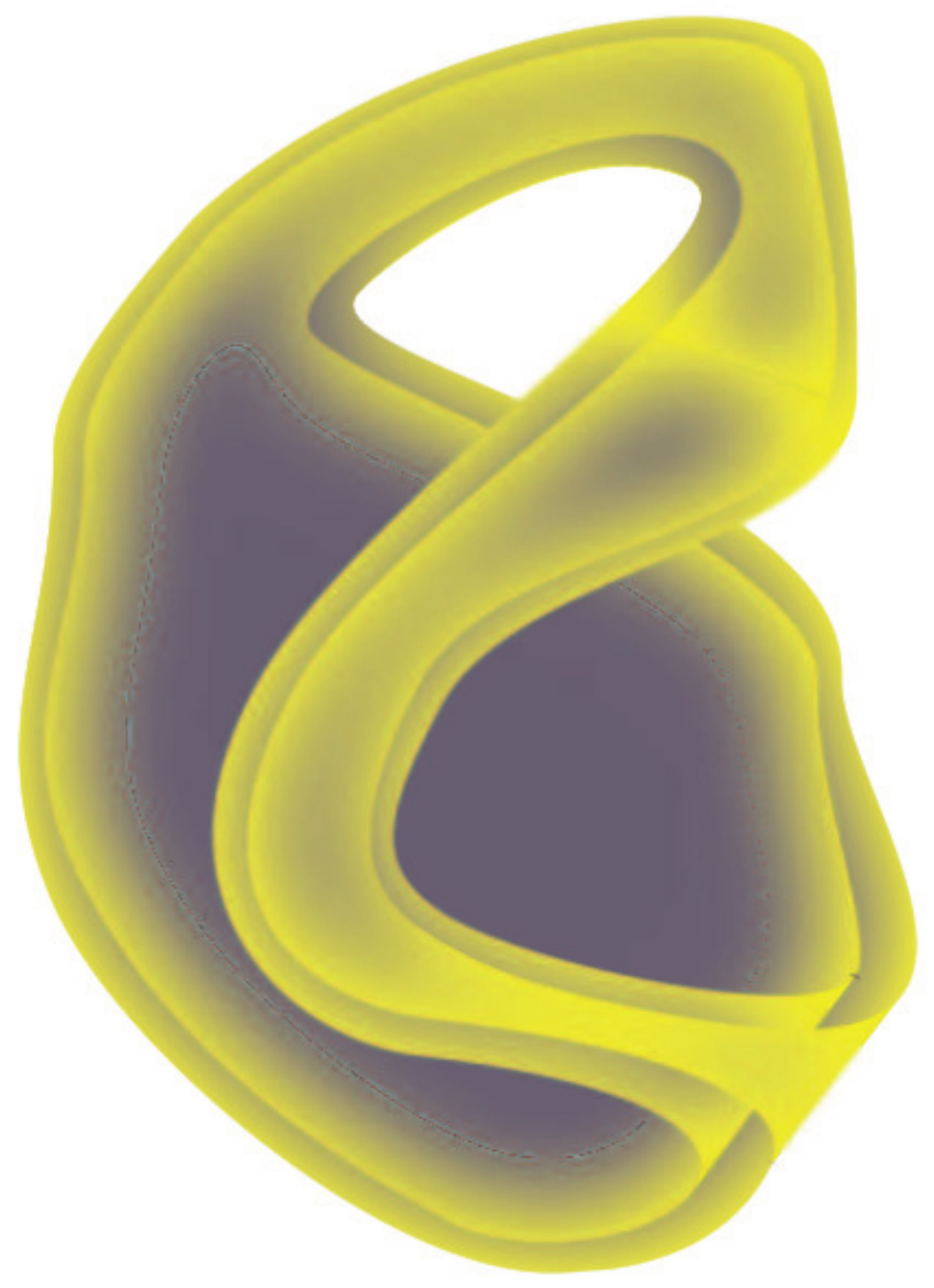} 
    \caption{
    (a)   A disc with a single clasp intersection (b) with right edge folded over and curled, revealing a 
    ($-4$)-half-twisted boundary annulus  (c), (d) The same annulus on a thickened M\"obius band as 2-fold cover. This is a twisted $I$-bundle in $\R^3$: a solid donut.}
    \label{fig:an-14}
\end{figure}
\begin{figure}[htbp] 
    \centering
(a)  \includegraphics[width=2.1in]{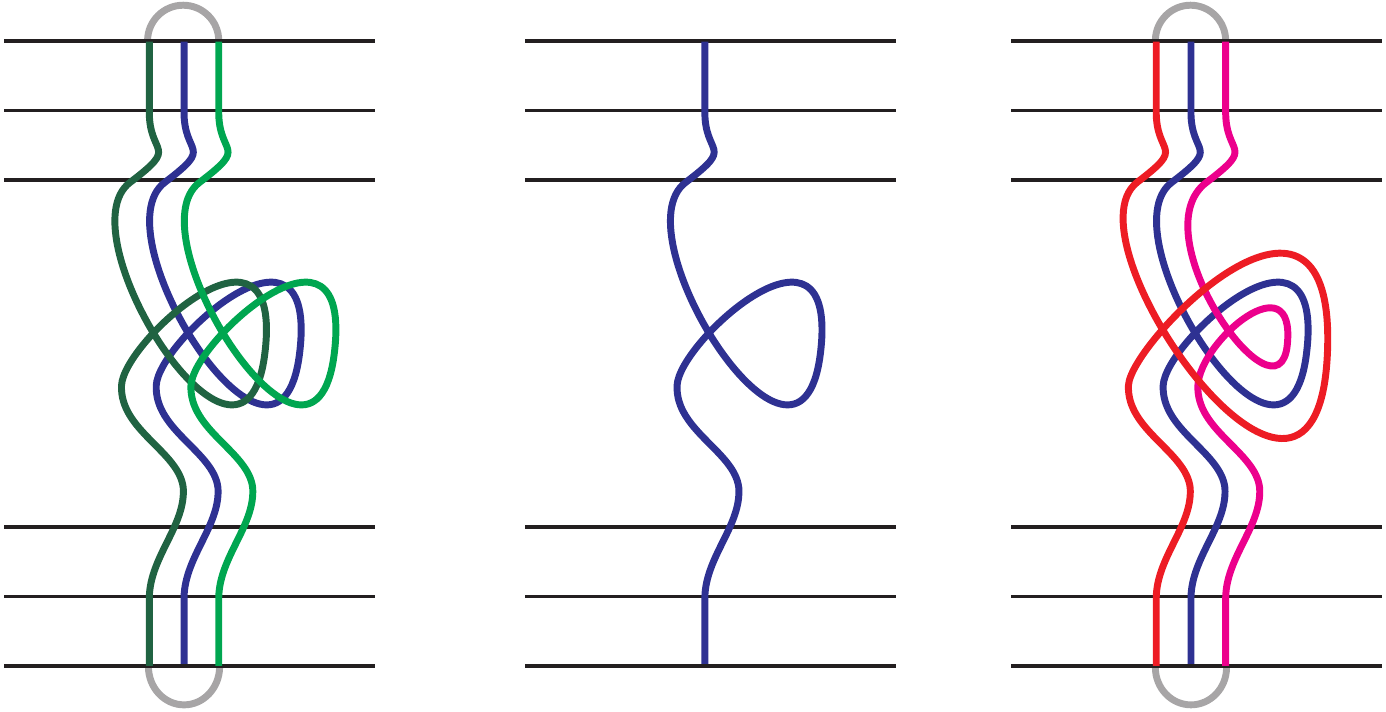} \qquad   \qquad (b) 
  \includegraphics[width=2.1in]{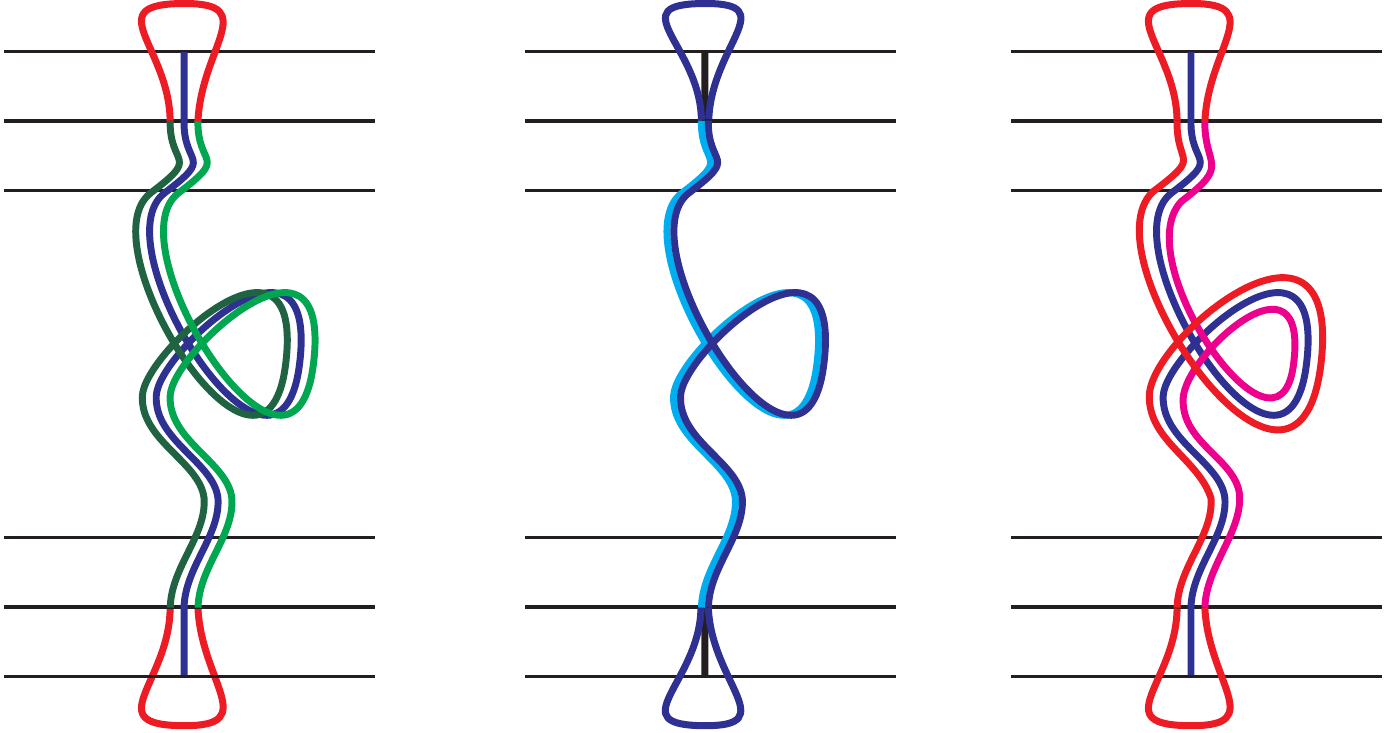} 
\vskip0.5cm
(c)    \includegraphics[width=4.5in]{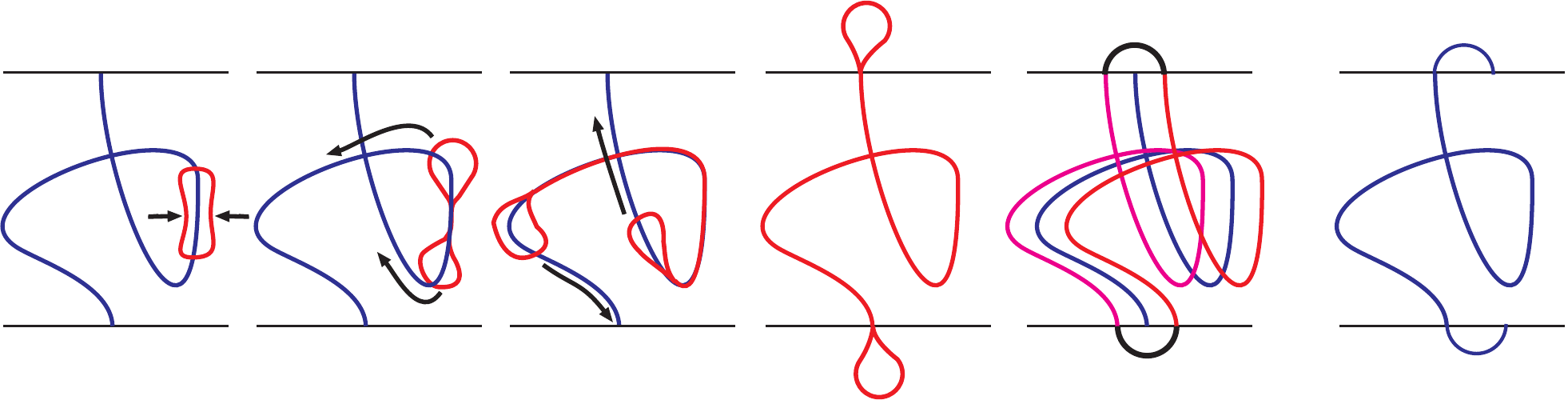} 

     \caption{
(a)  (Spin   pictures to obtain:/Cross sections of:)  An immersed disc (centre); an immersed sphere,   by 
  capping two rotated (translated) copies with an annulus (left);
  an immersed sphere as the boundary of an immersed regular neighbourhood (right).
  (b) Such spheres  are regularly homotopic. (c) An immersed disc has an immersed neighbourhood whose boundary sphere  is regularly homotopic to an embedding. The last picture indicates how to adjust an immersion so that it meets from the other side.}
    \label{fig:ArcsYes}
 \end{figure}

 \begin{figure}[htbp] 
    \centering
    \includegraphics[width= 5in]{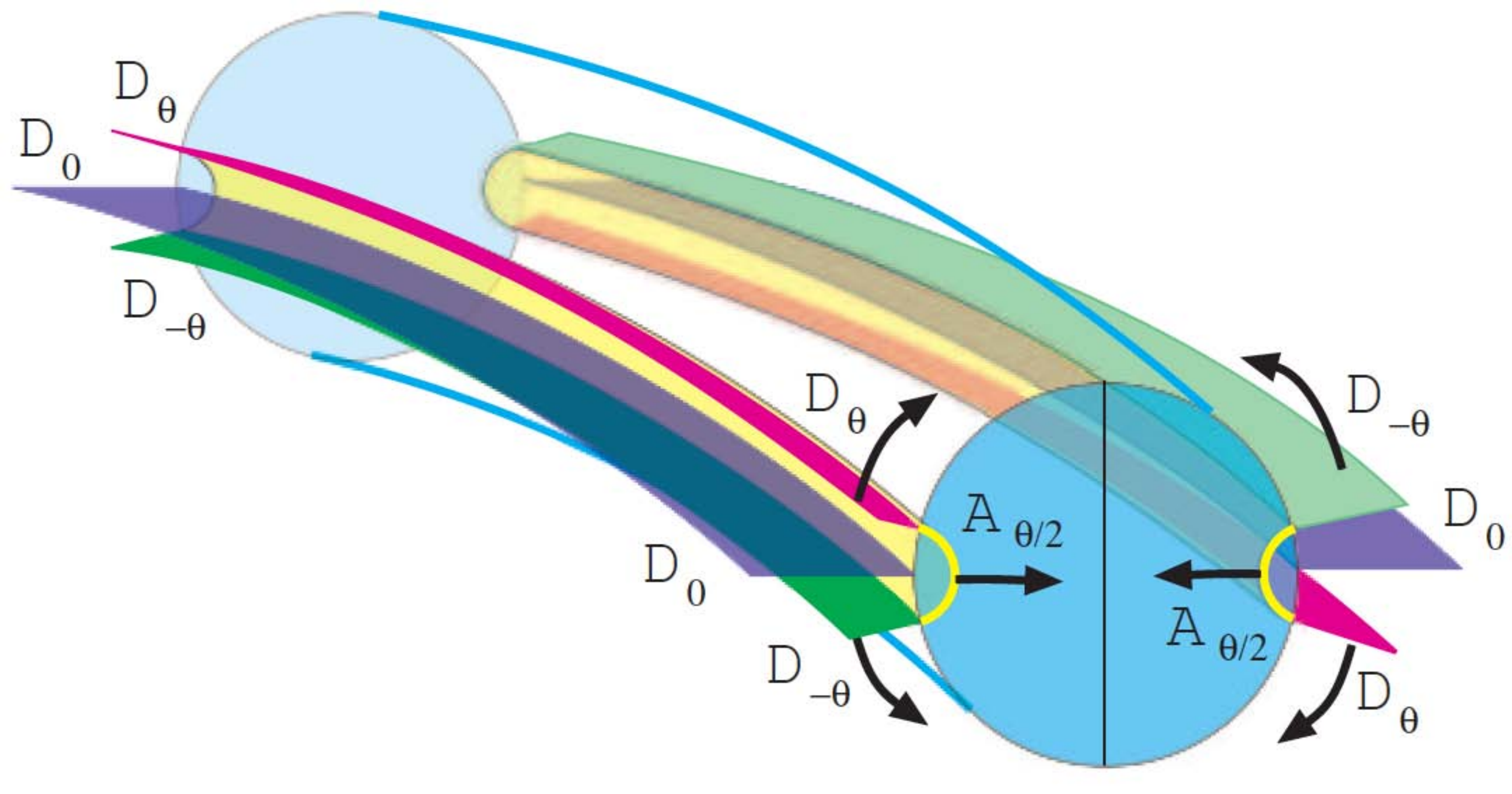} 
    \caption{ 
    Immersed spheres $S_\theta$ consist of two immersed discs $D_{\theta},\ D_{- \theta} \equiv D_{2\pi - \theta}$ and an annulus $A_{\theta/2}$.  Each annulus $A_{\theta/2}$ wraps twice around the solid donut. The two discs and two sheets of the annulus can be interchanged by increasing $\theta$, rotating the discs in opposite directions: the annulus passes smoothly through itself as it double covers the M\"obius band $M$. Discs inherit opposite orientations from $S_\theta$.}
    \label{fig: PushAnn1}
 \end{figure}
   
   \begin{figure}[htbp]  
       \centering
 \includegraphics[width=4in]{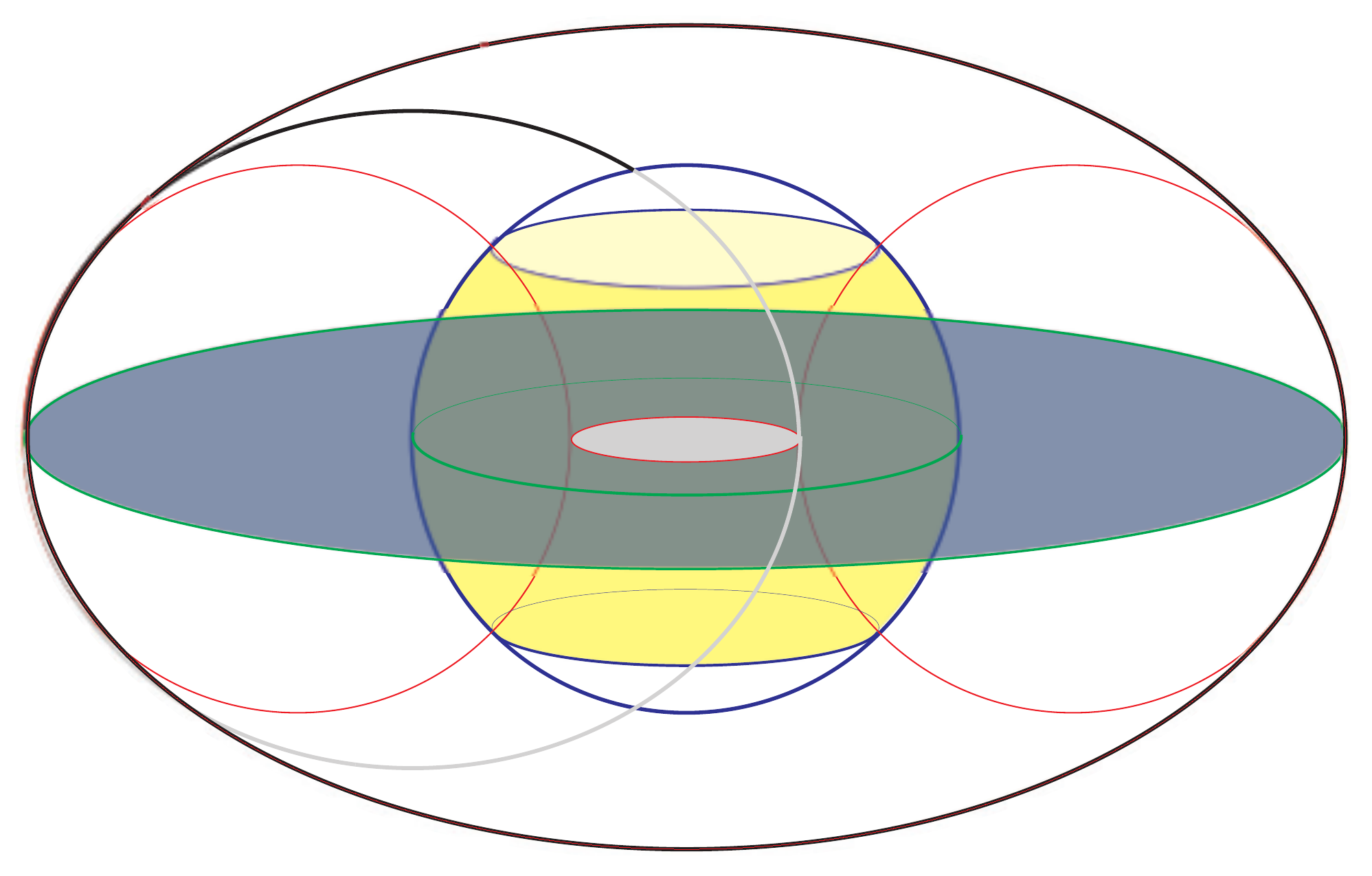} 
     \caption{ 
  After stereographic projection:   View in 3-dimensional space $\R^3$ of  the unit sphere $S$, Clifford torus $CT$, Hopf circle  and solid donut $SD$ with boundary torus $CT$ in $\R^3$: $A=S\cap SD$ and $SD\cap xy$-plane are  annuli meeting 
 along the unit circle. $D_0$ is the innermost $xy$-planar disc. The 3-ball and solid donut arise by rotating three coplanar unit discs, whose boundary circles meet orthogonally, around the $z$-axis. 
 The interior of the solid donut $SD$ can be foliated by annuli parallel to $A$; 
 similarly, the   complementary solid donut of $SD$ in $S^3$ can be foliated by discs `parallel' to $D_0$, which are caps of spheres, with a `disc'  $D^*$ the exterior of the $xy$-plane centered at infinity. These foliating discs and annuli can be chosen naturally to meet  $CT$   orthogonally, with boundary circles parallel to $\partial D_0 := \lambda$.  
    The boundary circle of $D_0$ can be parametrized in the standard way by angle $ \theta$; each point determines a unique Hopf circle, with $H$ as drawn corresponding to angle $\theta=0$.}
    \label{fig:SDCTS}
 \end{figure}

   \begin{figure}[htbp]  
    \centering
 \includegraphics[width=1.5in]{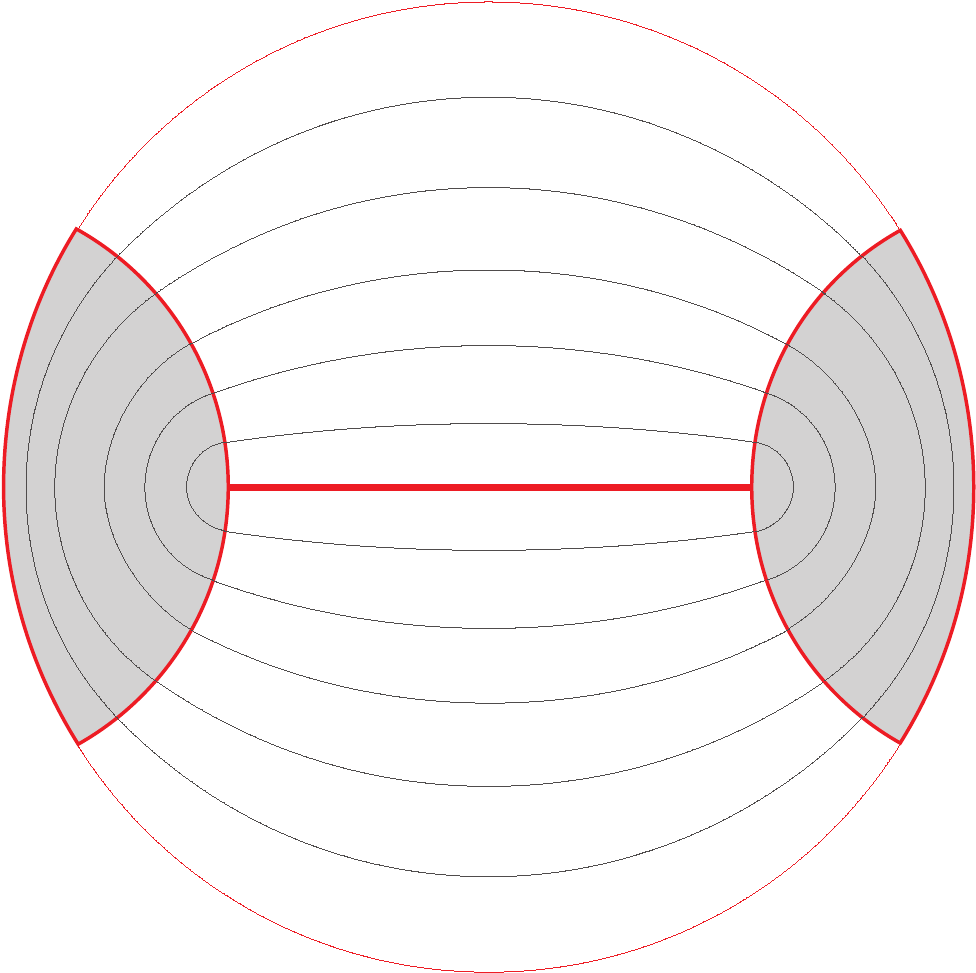} 
\qquad\qquad
 \includegraphics[width=1.5in]{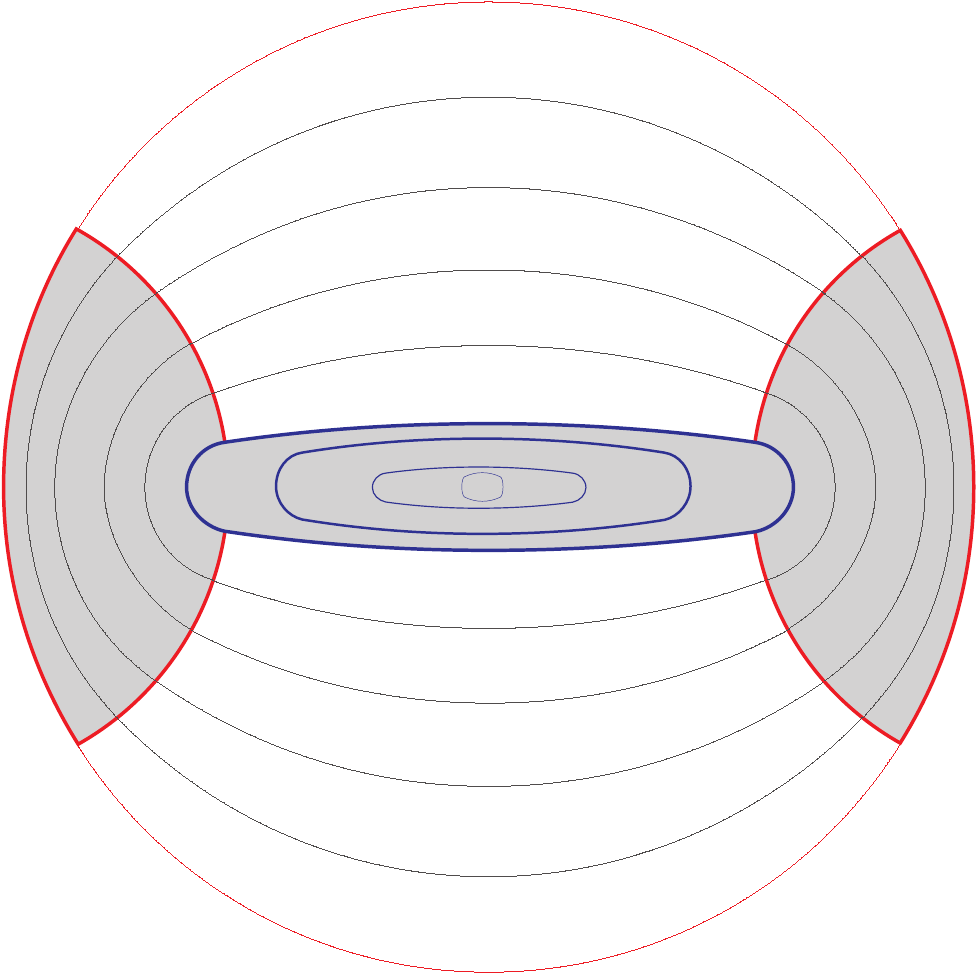} \qquad
  \includegraphics[width=1.5in]{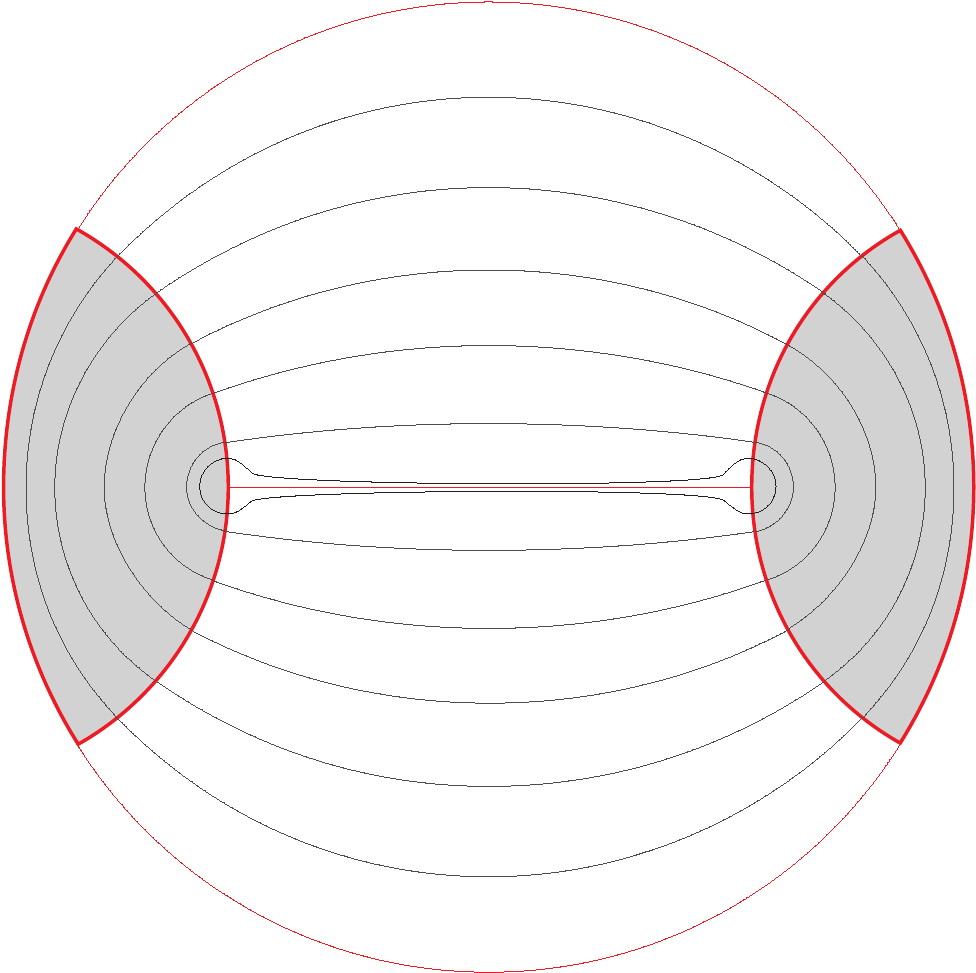} 
     \caption{ 
 Cross section of the unit ball intersecting the Clifford torus $CT$, $SD$ and $SD^*$, with spanning disc $D_0$ appearing as an interval: spin this around the vertical axis. Concentric spheres can shrink down to $D_0$, to a point, or to a shrinking neighbourhoods of a shrinking disc in $D_0$.}
    \label{fig:SphereExp}
 \end{figure}

   \begin{figure}[htbp]  
    \centering
 \includegraphics[width=5.2in]{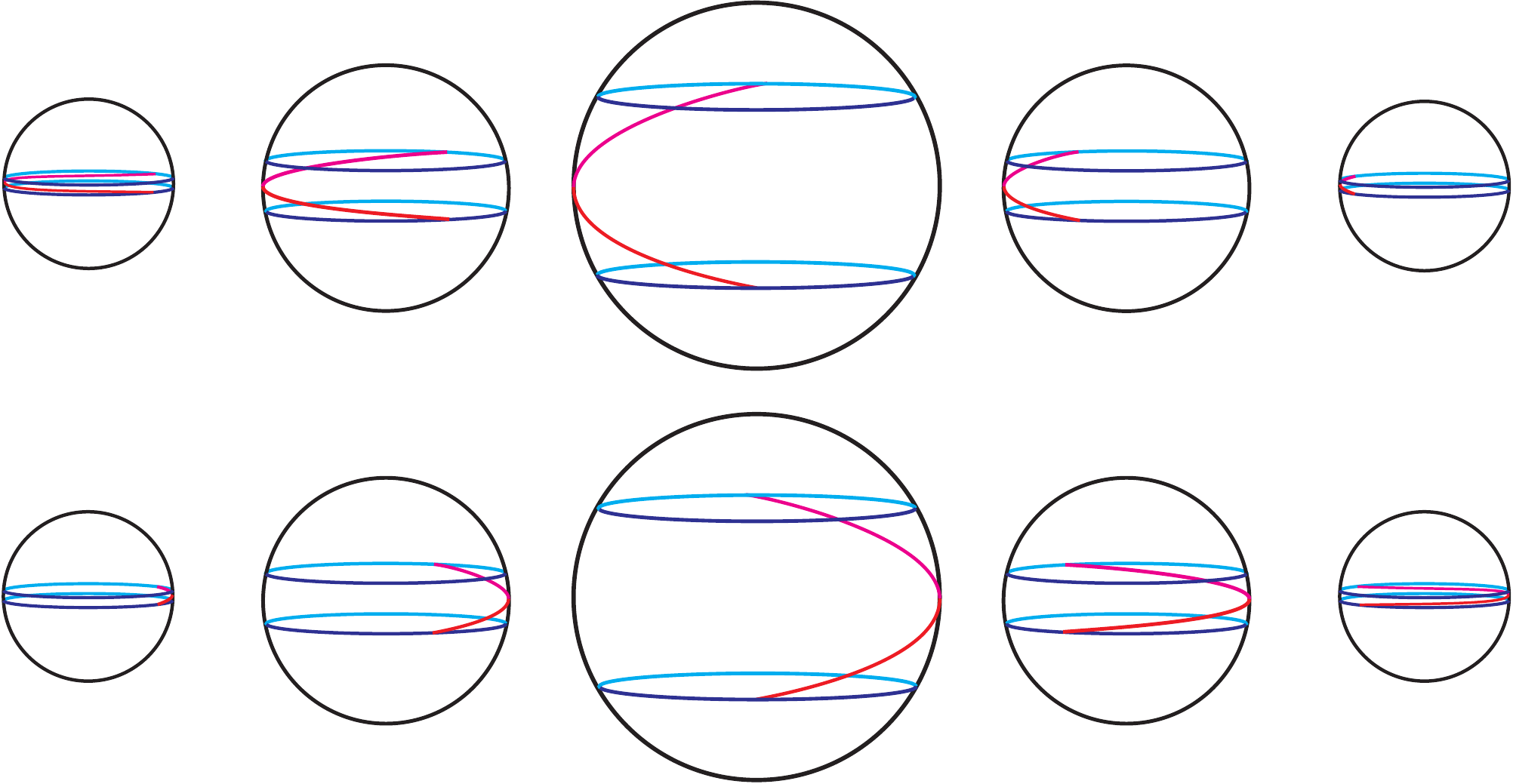} 
     \caption{ 
   Push $D_0$ around $SD^*$ in both directions, guided by the Hopf flow. The upper and lower spherical caps rotate in opposite directions. This pair of discs glues to an annulus in $SD$ which undergoes a complete Dehn twist. The top picture is the antipodal image in reverse order to the bottom, which, read left to right, shows the sequence of spheres beginning with a neighbourhood of $D_0$. }
    \label{fig:AntiTwist}
 \end{figure}

   \begin{figure}[htbp]
    \centering
 \includegraphics[width=6.2in]{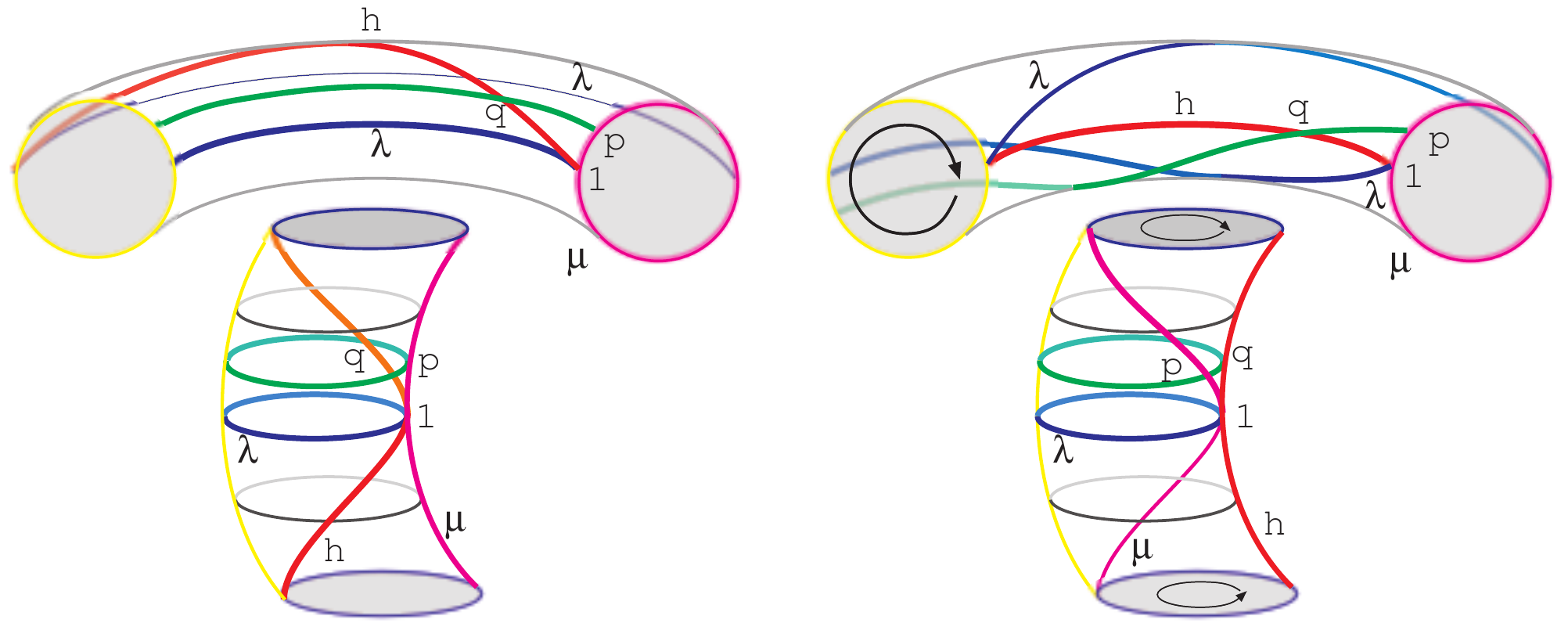} 
     \caption{
The left picture shows the standard longitude-meridian pair parametrizing the solid donuts $SD,\ SD^*$ as products $S^1\times D^2$; the right picture reveals the parametrization by the Hopf flow. 
 Only half of $SD$ is shown; $SD^*$ is shown cut along $D_\infty$. The point $P$ is indicated as `1'.}
    \label{fig:HopfFlowS}
 \end{figure}

   \begin{figure}[htbp] 
    \centering
 \includegraphics[width=6.2in]{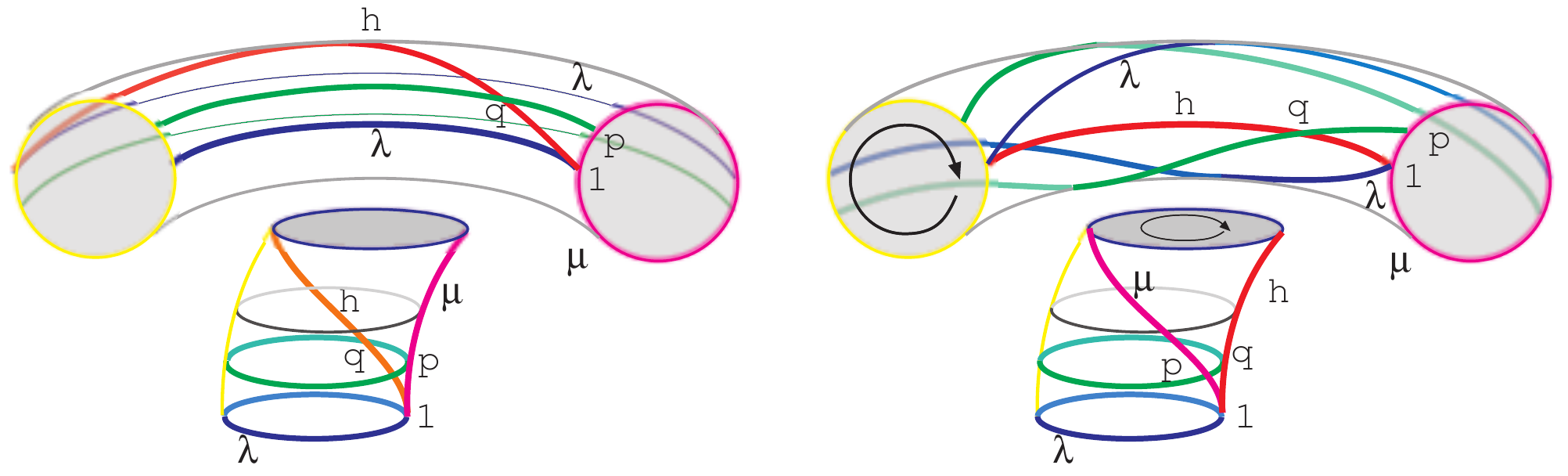} 
     \caption{
  The left picture shows the standard longitude-meridian pair parametrizing the solid donuts $SD/a,\ SD^*/a$ in $RP^3$, as products $S^1\times D^2$; with $\pi$-twist required for identification.  the right picture reveals the parametrization by the Hopf flow, enabling mapping $\R P^3$ into $\R^3$ so that the Hopf flow becomes rotation around the $z$-axis.   }
    \label{fig:HopfFlowRP}
 \end{figure}

\pagebreak
 	
\noindent
{Acknowledgement.} The author would like to express thanks to Professor Makoto Sakuma and Professor George Francis for their very helpful suggestions, and members of 
 the Department of Mathematics at Hiroshima University for their   hospitality.

E-mail: iain@ms.unimelb.edu.au
\end{document}